\documentclass[11pt,twoside]{article}%
\usepackage{enumerate}
\usepackage{enumitem}
\usepackage{amsmath,amssymb,amsfonts}
\usepackage{amsfonts}%
\usepackage{graphicx}
\usepackage{color}
\usepackage{pdfsync}
\providecommand{\U}[1]{\protect\rule{.1in}{.1in}}
\newtheorem{thm}{Theorem}[section]
\newtheorem{lem}[thm]{Lemma}

\newtheorem{prop}[thm]{Proposition}
\newtheorem{corol}[thm]{Corollary}
\newtheorem{claim}[thm]{Claim}

\newenvironment{pf}[1][\bfseries Proof]{\noindent{#1.} }{\hfill \rule{0.5em}{0.5em}\\}
\newcommand{\fim}{\hfill\rule{2mm}{2mm}}
\numberwithin{equation}{section}
\begin{document}
\title{Multiplicity and concentration of positive solutions for a class of quasilinear problems through Orlicz-Sobolev space}

\author{Claudianor O. Alves\thanks{C.O. Alves was partially supported by CNPq/Brazil  304036/2013-7  and INCT-MAT, \hspace*{.7cm} e-mail:  coalves@dme.ufcg.edu.br},~~  Ailton R. Silva\thanks{A.R. Silva,~ e-mail: ailton@dme.ufcg.edu.br}\\
Universidade Federal de Campina Grande\\
 Unidade Acad\^emica de Matem\'atica - UAMat\\
 CEP: 58.429-900 - Campina Grande - PB - Brazil}
\date{}
\maketitle

\begin{abstract}
In this paper, we study existence, multiplicity and concentration of positive solutions for the following class of quasilinear problems
\[
- \Delta_{\Phi}u + V(\epsilon x)\phi(\vert u\vert)u  =  f(u)\quad \mbox{in} \quad \mathbb{R}^{N} \,\,\, ( N\geq 2 ),
\]
where $\Phi(t) = \int_{0}^{\vert t\vert}\phi(s)sds$ is a N-function, $ \Delta_{\Phi}$ is the $\Phi$-Laplacian operator, $\epsilon$ is a positive parameter, $V : \mathbb{R}^{N} \rightarrow \mathbb{R} $ is a continuous function and $f : \mathbb{R} \rightarrow \mathbb{R} $ is a $C^{1}$-function.

\end{abstract}



\section{Introduction}

In this paper, we are concerned with existence, multiplicity and concentration of positive solutions for the following class of quasilinear problems
\begin{align}
\left\{
\begin{array}
[c]{rcl}%
- \Delta_{\Phi}u + V(\epsilon x)\phi(\vert u\vert)u & = & f(u)~  \mbox{in}~ \mathbb{R}^{N},\\
u \in W^{1, \Phi}(\mathbb{R}^{N}), &  &
\end{array}
\right. \tag{$ P_{\epsilon} $}\label{P1}%
\end{align}
where $ N\geq 2$, $\epsilon$  is a positive parameter, the operator $ \Delta_{\Phi} u  = \mathrm{div}(\phi(\vert\nabla u\vert)\nabla u) $,  where $\Phi(t) = \int_{0}^{\vert t\vert}\phi(s)sds$, named $\Phi$-Laplacian, is a natural extension of the $p$-Laplace operator, with $ p $ being a positive constant.

This class of problems arises in a lot of applications, such as \\

\noindent {\it Nonlinear Elasticity:} \, $ \Phi(t) = (1+t^{2})^{\alpha}-1, \alpha \in (1, \frac{N}{N-2})$, \\
\noindent { \it Plasticity:} \, $ \Phi(t) = t^{p}\ln(1+t), 1< \frac{-1+\sqrt{1+4N}}{2}<p<N-1, N\geq 3$, \\
\noindent { \it Generalized Newtonian Fluid:} \, $\Phi(t) = \int_{0}^{t}s^{1-\alpha}(\sinh^{-1}s)^{\beta}ds, 0\leq \alpha \leq 1$ and $\beta > 0$, \\
\noindent {\it Non-Newtonian Fluid:  } \, $\Phi(t) = \frac{1}{p}|t|^{p}$ for $p>1$,\\
\noindent { \it Plasma Physics: } \, $\Phi(t) = \frac{1}{p}|t|^{p} + \frac{1}{q}|t|^{q}$ where $1<p<q<N$ with $q \in (p, p^{*}).$ \\

\noindent The reader can find more details involving this subject in \cite{Db}, \cite{F1}, \cite{FN2} and their references.

In the last years,  problem (\ref{P1}) has received a special attention for the case where $\Phi(t)=\frac{1}{2}|t|^{2}$, that is, when it is of the form
$$
-\Delta{u}+V(\epsilon x)u=f(u) \,\,\, \mbox{in} \,\,\, \mathbb{R}^{N},
$$
where $f$ is a $C^1$-function with subcritical or critical growth and $V$ is a continuous function verifying some conditions. We believe that the two most famous conditions  on $V$ are the following: \\

\noindent {\bf Condition ${(V_0)}$} ( see  Rabinowitz \cite{rabinowitz} ): \\
$$
V_{\infty}:= \liminf_{\vert x\vert\rightarrow +\infty}V(x) > V_{0}:=\inf_{\mathbb{R}^{N}}V(x) >0.
$$

\noindent { \bf Condition $(V_1)$ }( see del Pino and Felmer \cite{DF1} and Gui \cite{G} ):  \\
$$
 V(x)\geq V_0=\inf_{x\in \mathbb{R}^N}V(x)>0 \ \ \hbox{for all}\ \ x\in \mathbb{R}^{N}, \eqno{(i)}
$$
and there is an open and bounded domain $\Lambda \subset \mathbb{R}^N$, such that
$$
  V_0<\inf_{x\in \partial\Lambda}V(x).\eqno{(ii)}
$$

In the literature, the reader can find a lot of papers, where the existence,  multiplicity and concentration of positive solutions have been obtained in connection with the geometry of the function $V$. In some of them, it was proved that the maximum points of the solutions are close to the set
$$
M = \{x \in \mathbb{R}^{N} \ : \ V(x) = V_{0} \},  \leqno{(M)}
$$
when $\epsilon$ is small enough. Moreover, in a lot of papers, the multiplicity of solution is also related to topology richness of the set $M$. The reader can find more information about this subject in \cite{AM}, \cite{DF1}, \cite{FW}, \cite{G}, \cite{O1},
 \cite{O2}, \cite{rabinowitz}  and their references.

For the case where $\Phi(t)=\frac{1}{p}|t|^{p}$ with $p>1$ or $\Phi(t)=\frac{1}{p}|t|^{p}+\frac{1}{p}|t|^{p}$ with $q,p>1$, the problem (\ref{P1}) becomes
$$
-\Delta_p{u}+V(\epsilon x)|u|^{p-2}u=f(u), \,\,\, \mbox{in} \,\,\, \mathbb{R}^{N}
$$
and
$$
-\Delta_p{u}-\Delta_q{u}+V(\epsilon x)(|u|^{p-2}u + |u|^{q-2}u)=f(u), \,\,\, \mbox{in} \,\,\, \mathbb{R}^{N}
$$
respectively. For these class of problems, the existence,  multiplicity and concentration of solutions have also been considered in some papers, see for example, \cite{AG2006}, \cite{AG1}, \cite{BB}, \cite{CV}, \cite{CEM}, \cite{F1}, \cite{LL} and their references.

For the case where the function $\Phi$ belongs to a large class of functions,  including for example,
$$
\Phi(t) = (1+t^{2})^{\alpha}-1 \,\,\,  \mbox{or} \,\,\, \Phi(t) = t^{p}\ln(1+t),
$$
we do not have found in the literature any paper where the existence,  multiplicity and concentration of solutions have been studied for some of the above functions, which appear in important applications as mentioned before. Motivated by this fact, in the present paper we study the existence,  multiplicity and concentration for a large class of functions $\Phi$ by supposing condition $(V_0)$. However, we would like point out that we are finishing other paper, where the same type of study has been made, by assuming that $V$ verifies $(V_1)$.

Actually, we have observed that there are interesting papers studying the existence of solution for (\ref{P1}) when $\epsilon=1$, we would like to cite the papers \cite{AGJ}, \cite{BBR2}, \cite{FN}, \cite{FN1}, \cite{FN2}, \cite{MR1}, \cite{MR2},  \cite{J} and references therein.

In this work, we generalize the arguments used in \cite{{AG2006}, AG1}, in the sense that, we will obtain the same type of results for a large class of operators. In the proof our result, we will work with N-function theory and Orlicz-Sobolev spaces. In the above references, the $L^{\infty}$-estimate of some sequences was obtained by using interaction Moser techniques, which does not work well when we are working with a general class of quasilinear operators. Here, we overcome this difficulty showing new estimates by adapting some arguments found in \cite{AF}, \cite{FUCHS}, \cite{Fusco}, \cite{LU} and \cite{ZF}.

Related to functions $\phi$ and $f$, our hypotheses are the following:
\begin{flushleft}
\textbf{Conditions on $\phi$:}
\end{flushleft}
The function $\phi: [0, +\infty)\rightarrow [0, +\infty)$ is a $C^{1}$- function satisfying
\begin{enumerate}[label={($\phi_\arabic{*}$})]
\setcounter{enumi}{0}
\item\label{H1} $\phi(t)$, $(\phi(t)t)^{'}>0$ for all $t>0$.

\item\label{H2} There exist $l,m \in (1, N)$ such that
  \[
  l\leq m <l^{*}=\frac{Nl}{N-l}
  \]
  and
  \[
  l\leq \frac{\phi(t)t^{2}}{\Phi(t)}\leq m \,\,\, \forall t\neq 0,
  \]
  where
  \[
  \Phi(t) = \int_{0}^{|t|} \phi(s)sds.
  \]
\item\label{H3} The function $\displaystyle\frac{\phi(t)}{t^{m-2}}$ is nonincreasing in $(0, +\infty)$.

\item\label{H4} The function $\phi$ is monotone.

\item\label{H5} There exists a constant $c>0$ such that
\[
\vert \phi^{'}(t)t\vert \leq c\phi(t) \quad \forall \ t \in [0, +\infty).
\]
\end{enumerate}

Hereafter, we will say that $\Phi \in \mathcal{C}_m$ if
$$
\Phi(t) \geq |t|^{m} \,\,\,\,\,\, \forall t \in \mathbb{R}. \leqno{(\mathcal{C}_m)}
$$
Moreover, let us denote by $\gamma$ the following real number
$$
\gamma=
\left\{
\begin{array}{l}
m, \,\,\,\, \mbox{if} \,\,\, \Phi \in \mathcal{C}_m, \\
\mbox{}\\
l, \,\,\,\, \mbox{if} \,\,\, \Phi \notin \mathcal{C}_m.
\end{array}
\right.
$$

Here, we would like to detach that the functions $\phi$ associated with each N-function mentioned in this introduction, except the N-function \linebreak $\Phi(t) = \int_{0}^{t}s^{1-\alpha}(\sinh^{-1}s)^{\beta}ds$, verify the conditions $(\phi_1)$-$(\phi_5)$.

\begin{flushleft}
\textbf{Conditions on $f$:}
\end{flushleft}
The function $f: \mathbb{R}\rightarrow \mathbb{R}$ is a $C^{1}$- function verifying:
\begin{enumerate}[label={($f_\arabic{*}$})]
  \item\label{f1}  There are functions $r,b:[0,+\infty) \to [0,+\infty)$ such that
  \[
  \limsup_{|t|\rightarrow 0} \frac{f^{'}(t)}{(r(|t|)|t|)^{'}}=0 \quad \mbox{and} \quad \limsup_{|t|\rightarrow +\infty} \frac{|f^{'}(t)|}{(b(|t|)|t|)^{'}}<+\infty.
  \]
  \item\label{f2} There exists $\theta > m$ such that
  \[
  0< \theta F(t) \leq f(t)t \quad \forall t>0,
  \]
  where
  \[
  F(t)= \int_{0}^{s}f(s)ds.
  \]
  \item\label{f3}  The function $\displaystyle\frac{f(t)}{t^{m -1}}$ is increasing for $t>0$.
\end{enumerate}
Related to functions $r $ and $b$ above, we assume that they are $C^{1}$- function  satisfying the ensuing conditions:
\begin{enumerate}[label={($r_\arabic{*}$})]
\setcounter{enumi}{0}
\item\label{R1} $r$ is increasing.
\item\label{R2} There exists a constant $\overline{c}>0$ such that
\[
\vert r^{'}(t)t \vert \leq \overline{c} r(t), \quad \forall \ t\geq 0.
\]
\item\label{R3} There exist positive constants $r_{1}$ and $r_{2}$ such that
 \[
  r_{1}\leq \frac{r(t)t^{2}}{R(t)}\leq r_{2}, \quad \forall t\neq 0,
  \]
  where
  \[
  R(t) = \int_{0}^{|t|} r(s)sds.
  \]
\item\label{R4} The function $R$ satisfies
\[
  \limsup_{t\rightarrow 0}\frac{R(t)}{\Phi(t)}< +\infty \quad \mbox{and} \quad \limsup_{|t|\rightarrow +\infty}\frac{R(t)}{\Phi_{*}(t)}=0.
  \]

\end{enumerate}

\begin{enumerate}[label={($b_\arabic{*}$})]
\setcounter{enumi}{0}
\item\label{B1} $b$ is increasing.
 \item\label{B2} There exists a constant $\widetilde{c}>0$ such that
\[
\vert b^{'}(t)t\vert \leq \widetilde{c}b(t), \quad \forall t\geq 0.
\]
\item\label{B3} There exist positive constants $b_{1}, b_{2} \in (1, \gamma^*)$ verifying
 \[
  b_{1}\leq \frac{b(t)t^{2}}{B(t)}\leq b_{2} \,\,\, \forall t\neq 0,
  \]
  where $\gamma^*=\frac{N\gamma}{N-\gamma}$ and
  \[
  B(t) = \int_{0}^{|t|} b(s)sds.
  \]
\item\label{B4} The function $B$ satisfies
\[
  \limsup_{t\rightarrow 0}\frac{B(t)}{\Phi(t)}< +\infty \quad \mbox{and} \quad \limsup_{|t|\rightarrow +\infty}\frac{B(t)}{\Phi_{*}(t)}=0,
  \]
\end{enumerate}
where $\Phi_{*}$ is the Sobolev conjugate function, which is defined by inverse function of
  \[
  G_{\Phi}(t) = \int_{0}^{t}\frac{\Phi^{-1}(s)}{s^{1+\frac{1}{N}}}ds.
  \]

Using the above hypotheses, we are able to state our main result.

\begin{thm}\label{T1}
Suppose that  \ref{H1}-\ref{H5}, \ref{R1}-\ref{R4}, \ref{B1}-\ref{B4}, \ref{f1}-\ref{f3} and $(V_0)$ hold. Then,
for any $\delta > 0$ small enough, there exists  $\epsilon_{\delta}> 0$ such that $(P_{\epsilon} )$ has at least $cat_{M_{\delta}}(M)$ positive solutions, for any $0 < \epsilon < \epsilon_{\delta}$, where
\[
M_{\delta} = \{x \in \mathbb{R}^{N} \ : \ dist(x,M) \leq \delta \}.
\]
Moreover, if $u_{\epsilon}$ denotes one of these solutions and $x_{\epsilon} \in \mathbb{R}^{N}$ is a global maximum point of $u_{\epsilon}$, then
\[
\lim_{\epsilon \rightarrow 0}V(\epsilon x_{\epsilon}) = V_{0}.
\]
\end{thm}

We would like point out that, if $Y$ is a closed subset of a topological space $X$, the Lusternik-Schnirelman category $cat_{X}(Y)$ is the least number of closed and contractible sets in $X$ which cover $Y$.

\section{Preliminaries about Orlicz-Sobolev spaces}

\label{expoents_variaveis} In this section, we recall some properties of Orlicz and Orlicz-Sobolev spaces. We refer to \cite{Adams, FN, Rao} for the fundamental properties of these spaces.
First of all, we recall that a continuous function $\Phi : \mathbb{R} \rightarrow [0,+\infty)$ is a
N-function if:
\begin{description}
  \item[$(i)$] $\Phi$ is convex.
  \item[$(ii)$] $\Phi(t) = 0 \Leftrightarrow t = 0 $.
  \item[$(iii)$] $\displaystyle\lim_{t\rightarrow0}\frac{\Phi(t)}{t}=0$ and $\displaystyle\lim_{t\rightarrow+\infty}\frac{\Phi(t)}{t}= +\infty$ .
  \item[$(iv)$] $\Phi$ is even.
\end{description}
We say that a N-function $\Phi$ verifies the $\Delta_{2}$-condition, denote by  $\Phi \in \Delta_{2}$, if
\[
\Phi(2t) \leq K\Phi(t)\quad \forall t\geq 0,
\]
for some constant $K > 0$. In what follows, fixed an open set $\Omega \subset \mathbb{R}^{N}$ and a N-function $\Phi$, we define the Orlicz space associated with $\Phi$ as
\[
L^{\Phi}(\Omega) = \left\{  u \in L_{loc}^{1}(\Omega) \colon \ \int_{\Omega} \Phi\Big(\frac{|u|}{\lambda}\Big)dx < + \infty \ \ \mbox{for some}\ \ \lambda >0 \right\}.
\]
The space $L^{\Phi}(\Omega)$ is a Banach space endowed with Luxemburg norm given by
\[
\Vert u \Vert_{\Phi} = \inf\left\{  \lambda > 0 : \int_{\Omega}\Phi\Big(\frac{|u|}{\lambda}\Big)dx \leq1\right\}.
\]
The complementary function $\widetilde{\Phi}$ associated with $\Phi$ is given
by the Legendre's transformation, that is,
\[
\widetilde{\Phi}(s) = \max_{t\geq 0}\{ st - \Phi(t)\} \quad  \mbox{for} \quad s\geq0.
\]
The functions $\Phi$ and $\widetilde{\Phi}$ are complementary each other. In \cite{FN, Rao}, we find the ensuing property involving   $\Phi, \widetilde{\Phi}$:
$$
\Phi, \widetilde{\Phi} \in \Delta_{2} \,\,\, \mbox{if, and only if,} \,\, (\phi_2) \,\,\, \mbox{holds}.
$$
Moreover, we also have a Young type inequality given by
\[
st \leq \Phi(t) + \widetilde{\Phi}(s)\quad \forall s, t\geq0.
\]
Using the above inequality, it is possible to prove a H\"older type inequality, that is,
\[
\Big| \int_{\Omega}uvdx \Big| \leq 2 \Vert u \Vert_{\Phi}\Vert v \Vert_{\widetilde{\Phi}}\quad \forall u \in L^{\Phi}(\Omega) \quad \mbox{and} \quad \forall v \in L^{\widetilde{\Phi}}(\Omega).
\]
The corresponding Orlicz-Sobolev space is defined as
\[
W^{1, \Phi}(\Omega) = \Big\{ u \in L^{\Phi}(\Omega) \ :\ \frac{\partial u}{\partial x_{i}} \in L^{\Phi}(\Omega), \quad i = 1, ..., N\Big\}.
\]
Usually, the Orlicz-Sobolev space is obtained by completion of $C_{0}^{\infty}(\Omega)$ with respect to norm
\[
\Vert u \Vert = \Vert \nabla u \Vert_{\Phi} + \Vert u \Vert_{\Phi}.
\]
The spaces $L^{\Phi}(\Omega)$ and $W^{1, \Phi}(\Omega)$ are separable and reflexive when $\Phi$ and $\widetilde{\Phi}$ satisfy the $\Delta_{2}$-condition. The condition $\Delta_{2}$ implies that
\[
u_{n}\to u \ \mbox{in} \ \ L^{\Phi}(\Omega)\quad \Leftrightarrow \quad \int_{\Omega}\Phi(\vert u_{n} - u\vert)dx \rightarrow 0
\]
and
\[
u_{n}\to u \ \mbox{in} \ \ W^{1, \Phi}(\Omega)\quad \Leftrightarrow  \quad \int_{\Omega}\Phi(\vert \nabla u_{n} - \nabla u\vert)dx \to 0\,\, \mbox{and}\,\, \int_{\Omega}\Phi(\vert u_{n} - u\vert)dx \to 0.
\]
In the literature, we have some important embeddings involving the Orlicz-Sobolev spaces. In \cite{Adams, DT}, it is showed that if $A$ is a N-function satisfying
\[
\limsup_{t\to 0}\frac{A(t)}{\Phi(t)}< +\infty \quad \mbox{and}\quad \limsup_{t\rightarrow +\infty}\frac{A(t)}{\Phi_{*}(t)}< +\infty,
\]
then the embedding
\[
W^{1, \Phi}(\Omega)\hookrightarrow L^{A}(\Omega)
\]
is continuous. If $\Omega$ is a bounded domain, the embedding is compact.

Before to conclude this introduction, we will prove a more general result than Theorem 1.2 found in \cite{ACG}. In our work, the hypothesis \ref{H5} is weaker that $(\phi_{3})^{'}$ and we do not use the hypotheses $(\phi_{2})^{'}$ from that paper. Here, we can work with the function $\phi(t) = 2\alpha(1+t^{2})^{\alpha-1}$, which cannot be considered in \cite{ACG}.
\begin{prop}\label{Teoape1}
Let $\Omega \subset \mathbb{R}^{N}$ be a smooth domain. Assume \ref{H1}-\ref{H5}. Let $(\eta_{n})$ be a sequence of vector fields, $\eta_{n} : \Omega \rightarrow \mathbb{R}^{N}$ satisfying
\[
\eta_{n} \in L^{\Phi}(\Omega)\times\ldots\times L^{\Phi}(\Omega)\quad \mbox{and}\quad\eta_{n}(x)\rightarrow 0 \ \mbox{a.e.} \ x \in \Omega.
\]
If $(\Vert \eta_{n}\Vert_{\Phi})$ is bounded and $a(y) := \phi(\vert y\vert)y$ for all $y \in \mathbb{R}^{N}$, then
\[
\int_{\Omega}\widetilde{\Phi}(\vert a(\eta_{n} + w) - a(\eta_{n}) - a(w)\vert)dx = o_{n}(1),
\]
for each $w \in L^{\Phi}(\Omega)\times\ldots\times L^{\Phi}(\Omega)$.
\end{prop}
\begin{pf} The proof is essentially the same explored in \cite[Theorem 2.1]{ACG}. If $\phi$ is decreasing, the proof is the same. However, if $\phi$ is non-decreasing, it is necessary to make a little adjustment. As done in \cite{ACG}, the key point in the proof is to show that, given $\varepsilon>0$, there exist $C, c_\varepsilon>0$ such that
\begin{equation}\label{VZ}
\vert a(\eta_{n} + w) - a(\eta_{n})\vert \leq C\big( \varepsilon \phi(\vert \eta_{n}(x)\vert)\vert \eta_{n}(x)\vert + c_\varepsilon\phi( \vert w(x) \vert)\vert w(x) \vert \big),
\end{equation}
where $C$ is independent of $\varepsilon$.

We can prove this fact using only the hypotheses \ref{H1}-\ref{H5}. In fact, just note that for each $w \in L^{\Phi}(\Omega)\times\ldots\times L^{\Phi}(\Omega)$, the assumptions \ref{H1}-\ref{H5} imply
\begin{eqnarray}\label{VZ1}
\vert a_{i}(\eta_{n} + w) - a_{i}(\eta_{n})\vert&\leq&C_1\big(\vert w \vert\phi(\vert\eta_{n} \vert) + \vert w \vert\phi(\vert w \vert)\big),
\end{eqnarray}
for some $C_1>0$. Define the sequence $H_{n} = C_1\vert w \vert \phi(\vert \eta_{n}\vert)$.
\begin{claim}\label{AF12}
For all $\varepsilon>0$, there exist $C_2, c_\varepsilon >0$, such that
\begin{eqnarray}\label{VZ2}
H_{n}(x) \leq C_2\big( \varepsilon \phi(\vert \eta_{n}(x)\vert)\vert \eta_{n}(x)\vert + c_\varepsilon\phi( \vert w(x) \vert)\vert w(x) \vert \big) \,\,\, \forall n \in \mathbb{N}.
\end{eqnarray}
\end{claim}
Indeed, fixed $x \in \mathbb{R}^{N}$, we set
$$
A_{x} = \{ n \in \mathbb{N} : \ \vert w(x) \vert \leq \vert \eta_{n}(x)\vert \} \,\,\, \mbox{and} \,\,\, B_{x}:= \{ n \in \mathbb{N} : \ \vert \eta_{n}(x)\vert  \leq \vert w(x)\vert\}.
$$
If $n \in A_{x}$, \cite[Lemma 2.2]{AGJ} and $\Delta_{2}$-condition lead to
\begin{eqnarray*}
H_{n}(x)\vert \eta_{n}(x)\vert &=& \frac{c\vert w(x) \vert}{\varepsilon}\big(\varepsilon \phi(\vert \eta_{n}(x)\vert)\vert \eta_{n}(x)\vert\big)\\
&\overbrace{\leq}^{Young}& \Phi(\frac{c}{\varepsilon}\vert w(x)\vert) + \widetilde{\Phi}(\varepsilon \phi(\vert \eta_{n}(x)\vert)\vert \eta_{n}(x)\vert)\\
&\leq& C_{\varepsilon}\Phi(\vert w(x)\vert) + \varepsilon \Phi(2\phi(\vert \eta_{n}(x)\vert))\\
&\leq& C_{\varepsilon}\Phi(\vert w(x)\vert) + C_3\varepsilon\Phi(\vert \eta_{n}(x)\vert)\\
&\leq& \frac{1}{l}\big( C_{\varepsilon}\phi( \vert w(x) \vert)\vert w(x) \vert^{2} + C_3\varepsilon \phi(\vert \eta_{n}(x)\vert)\vert \eta_{n}(x)\vert^{2}\big),
\end{eqnarray*}
which implies that there exists $C_4>0$ such that
\begin{eqnarray*}
H_{n}(x) \leq C_4\big( \varepsilon \phi(\vert \eta_{n}(x)\vert)\vert \eta_{n}(x)\vert + c_\varepsilon\phi( \vert w(x) \vert)\vert w(x) \vert \big).
\end{eqnarray*}
Now, if $n \in B_{x}$, using the fact that $\phi$ is increasing, we have that
$$
H_{n}(x) \leq C_2\phi( \vert w(x) \vert)|w(x)|,
$$
and so, there exists $C_5>0$ such that
\begin{eqnarray*}
H_{n}(x) \leq C_5\big( \varepsilon \phi(\vert \eta_{n}(x)\vert)\vert \eta_{n}(x)\vert + c_\varepsilon\phi( \vert w(x) \vert)\vert w(x) \vert \big),
\end{eqnarray*}
showing the Claim \ref{AF12}. By \eqref{VZ1} and \eqref{VZ2}, we get \eqref{VZ}.
\end{pf}

\section{Autonomous case}

In this section, our goal is to prove the existence of positive ground state solution for the problem
\begin{align}
\left\{
\begin{array}
[c]{rcl}%
- \Delta_{\Phi}u + \mu\phi(|u|)u & = & f(u)~  \mbox{in}~ \mathbb{R}^{N},\\
u \in W^{1, \Phi}(\mathbb{R}^{N}), &  &
\end{array}
\right. \tag{$ P_{\mu} $}\label{PA}%
\end{align}
where $\mu$ is a positive parameter and the functions $\phi$ and $f$ verify the hypotheses mentioned in the introduction. We recall that $u \in W^{1,\Phi}(\mathbb{R}^{N})$ is a solution ( weak solution ) of (\ref{PA}) if
$$
\int_{\mathbb{R}^{N}}\phi(|\nabla u|)\nabla u \nabla vdx + \mu \int_{\mathbb{R}^{N}}\phi(|u|)uvdx = \int_{\mathbb{R}^{N}}f(u)vdx \,\,\, \forall v \in W^{1,\Phi}(\mathbb{R}^{N}).
$$

Here, we would like point out that the above class of problems was considered in \cite{AGJ}. However, in that paper, it was not established the existence of ground state solution, characterization of the mountain pass level as the minimum of the energy function on the Nehari, regularity and positiveness of the solutions. This study is crucial for us, because some of our estimates depend these information.

Since we intend to find positive solutions, we will assume that
\begin{equation} \label{POSITIVIDADE}
f(t) = 0 \quad \mbox{for all} \ t< 0.
\end{equation}

From now on, $Y_{\mu}$ denotes the Orlicz-Sobolev space $W^{1, \Phi}(\mathbb{R}^{N})$ endowed with the norm
\[
\Vert u\Vert_{Y_{\mu}} = \Vert \nabla u\Vert_{\Phi} + \mu\Vert u\Vert_{\Phi}.
\]
It is easy to check that embeddings
\[
Y_{\mu}\hookrightarrow L^{\Phi}( \mathbb{R}^{N}) \quad \mbox{and}\quad Y_{\mu}\hookrightarrow L^{B}( \mathbb{R}^{N})
\]
are continuous. Moreover, let us denote by $E_{\mu} : Y_{\mu} \rightarrow \mathbb{R}$ the energy
functional related to \eqref{PA} given by
\[
E_{\mu}(u)= \int_{\mathbb{R}^{N}}\Phi(|\nabla u|)dx + \mu \int_{\mathbb{R}^{N}}\Phi(|u|)dx - \int_{\mathbb{R}^{N}}F(u)dx.
\]
Using standard arguments, it follows that $E_{\mu} \in C^{1}(Y_{\mu}, \mathbb{R})$ with
\[
E^{'}_{\mu}(u)v = \int_{\mathbb{R}^{N}}\phi(|\nabla u|)\nabla u \nabla vdx + \mu \int_{\mathbb{R}^{N}}\phi(|u|)uvdx - \int_{\mathbb{R}^{N}}f(u)vdx,
\]
for all $u,v \in Y_{\mu}$.

Note that $u \in Y_{\mu}$ is a weak solution of (\ref{PA}) if, and only if, $u$ is a critical point of $E_\mu$.  Moreover, it is very important to observe that if $u \in Y_{\mu}$ is a critical point of $E_\mu$, then $u$ is nonnegative,   because the equality $E'(u)u_{-} =0$ yields $u_{-}=0$, where $u_{-}=\max\{-u,0\}$.

Arguing as in \cite[Lemma 4.1]{AGJ}, the functional $E_{\mu}$ verifies the mountain pass geometry. Then, we can apply a version of  Mountain Pass Theorem without the Palais-Smale condition found in \cite{Willen}, which assures the existence of a $(PS )_{d_{\mu}}$
sequence $(u_{n}) \subset Y_{\mu}$, that is, a sequence satisfying
\[
E_{\mu}(u_{n})\rightarrow d_{\mu} \quad \mbox{and}\quad E_{\mu}^{'}(u_{n})\rightarrow 0,
\]
where $d_{\mu}$ is the mountain level of $E_{\mu}$, that is,
\[
d_{\mu}= \inf_{\alpha \in \Gamma}\max_{t \in [0, 1]}E_{\mu}(\alpha(t))
\]
with
\[
\Gamma := \{ \alpha \in C\big( [0, 1], Y_{\mu}\big);\, E_{\mu}(0) = 0 \ \ \mbox{and} \ \  E_{\mu}(\alpha(1))<0 \}.
\]

Next, we will show some characterizations of $d_{\mu}$, which will  be used later on. To this end, we denote by $\mathcal{M}_{\mu}$ the Nehari manifold associated to $E_{\mu}$ given by
\[
\mathcal{M}_{\mu} = \big\{ u \in Y_{\mu} \backslash \{0\} \ : \ \ E^{'}_{\mu}(u)u = 0 \big\}.
\]
Furthermore, we denote by $d_{\mu, 1}$ and $d_{\mu, 2}$ the ensuing real numbers
\[
d_{\mu, 1} = \inf_{u \in \mathcal{M}_{\mu}}E_{\mu}(u) \quad \mbox{and} \quad d_{\mu, 2} = \inf_{u \in Y_{\mu}\backslash\{0\}}\max_{t\geq 0}E_{\mu}(tu).
\]

\begin{lem}\label{Lcar}
Assume that \ref{H1}-\ref{H3}, \ref{R1}-\ref{R4}, \ref{B1}-\ref{B4} and \ref{f1}-\ref{f3} hold. Then, for each $u \in Y_{\mu}\backslash \{0\}$, there exists a unique $t_{u}>0$ such that $t_{u}u \in \mathcal{M}_{\mu}$ and $E_{\mu}(t_{u}u) = \displaystyle\max_{t \geq 0}E_{\mu}(tu)$. Moreover,
\[
d_{\mu} = d_{\mu, 1} = d_{\mu, 2}.
\]
\end{lem}
\begin{pf}
For each  $u \in Y_{\mu}\backslash \{0\}$, we define $g(t)= E_{\mu}(tu)$, that is,
\begin{eqnarray*}
g(t) = \displaystyle\int_{\mathbb{R}^{N}}\Phi(\vert\nabla (tu)\vert)dx + \mu\displaystyle\int_{\mathbb{R}^{N}}\Phi(\vert tu\vert)dx - \displaystyle\int_{\mathbb{R}^{N}}F(tu)dx.
\end{eqnarray*}
By a direct computation, it follows that $g(t)>0$ for $t$ sufficiently small and $g(t)<0$ for $t$ sufficiently large. Thereby, there is $t_{u}>0$ such that
\begin{eqnarray*}
g(t_{u}) = \max_{t\geq 0}g(t) = \max_{t \geq 0}E_{\mu}(tu),
\end{eqnarray*}
implying that $t_{u}u \in \mathcal{M}_{\mu}$.

Note that if $u\in \mathcal{M}_{\mu}$, then $u_{+}=\max\{u, 0\}\neq 0$. Suppose that there exist $t_{1}, t_{2} >0$ such that $t_{1}u, t_{2}u \in \mathcal{M}_{\mu}$. Then
\[
\displaystyle \int_{\mathbb{R}^{N}}\phi(|\nabla (t_{1}u)|)|\nabla(t_{1}u)|^{2}dx + \mu \displaystyle\int_{\mathbb{R}^{N}}\phi(|t_{1}u|)|t_{1}u|^{2}dx = \displaystyle\int_{[u>0]}f(t_{1}u)t_{1}udx
\]
and
\[
\displaystyle \int_{\mathbb{R}^{N}}\phi(|\nabla (t_{2}u)|)|\nabla(t_{2}u)|^{2}dx + \mu \displaystyle\int_{\mathbb{R}^{N}}\phi(|t_{2}u|)|t_{2}u|^{2}dx = \displaystyle\int_{[u>0]}f(t_{2}u)t_{2}udx.
\]
Thus,
\begin{eqnarray*}
 &&\displaystyle\int_{\mathbb{R}^{N}}\Bigg[\displaystyle\frac{\phi(t_{1}|\nabla u|)}{(t_{1}|\nabla u|)^{m-2}} - \displaystyle\frac{\phi(t_{2}|\nabla u|)}{(t_{2}|\nabla u|)^{m-2}}\Bigg]|\nabla u|^{m}dx\\
&+&\mu\displaystyle\int_{\mathbb{R}^{N}}\Bigg[\displaystyle\frac{\phi(t_{1}| u|)}{(t_{1}|u|)^{m-2}} - \displaystyle\frac{\phi(t_{2}|u|)}{(t_{2}|u|)^{m-2}}\Bigg]|u|^{m}dx = \displaystyle\int_{[u>0]}\Bigg[\displaystyle\frac{f(t_{1}u)}{(t_{1}u)^{m-1}}-\displaystyle\frac{f(t_{2}u)}{(t_{2}u)^{m-1}}\Bigg]u^{m}dx.
\end{eqnarray*}
From \ref{H3} and \ref{f3}, $t_{1} = t_{2}$. Now, the proof follows repeating the same arguments found in \cite[Theorem 4.2]{Willen}.
\end{pf}

The next lemma establishes some preliminary properties involving the Palais-Smale sequences of $E_\mu$.

\begin{lem}\label{la}
Let $(u_{n})$ be a $(PS)_{c}$ sequence for $E_\mu$ . Then, $(u_{n})$ is a bounded sequence in $Y_{\mu}$. Moreover, if $u_n \rightharpoonup u$ in $Y_{\mu}$, for some subsequence, we have that
\begin{description}
  \item[$(i)$] $ u_{n}(x)\rightarrow  u(x)\quad \mbox{a. e. in} \ \mathbb{R}^{N}$,
  \item[$(ii)$] $ \nabla u_{n}(x)\rightarrow \nabla u(x)\quad \mbox{a. e. in} \ \mathbb{R}^{N}$,
  \item[$(iii)$] $E^{'}_{\mu}(u) = 0$,
  \item[$(iv)$] ${u_{n}}_{-} = \max\{-u_{n}, 0\}\rightarrow 0$ in $Y_{\mu}$,
	\item[$(v)$] ${u_{n}}_{+} = \max\{u_{n}, 0\}$ is a $(PS)_c$ sequence for $E_\mu$. Hence, we can assume that $u_n \geq 0$ for all $n \in \mathbb{N}$, and  so, $u \geq 0$.
\end{description}
\end{lem}
\begin{pf}
The proof of $(i), (ii)$ and $(iii)$ follow using similar argument explored in \cite[Lemma 4.3]{AGJ}. The proof of $(iv)$ follows as in \cite[Lemma 2.3]{AG1}, while that $(v)$ is an immediate consequence of $(iv)$.
\end{pf}
\begin{lem}\label{Lemacompaut}
Let $(u_{n}) \subset Y_{\mu}$ be a $(PS)_{c}$ sequence for $E_{\mu}$ with $u_{n}\rightharpoonup 0 $ in $Y_{\mu}$. Then,
\begin{description}
 \item[$a)$] $u_{n}\rightarrow 0$ in $Y_{\mu}$,
\end{description}	
or
\begin{description}
  \item[$b)$] There exist constants $R_0, \tau>0$ such that
$$
\liminf_{n \rightarrow +\infty}\sup_{y \in \mathbb{R}^{N}}\int_{B_{R_0}(y)}\Phi(\vert u_{n}\vert)dx \geq \tau >0.
$$
\end{description}	

\end{lem}
\begin{pf}
Suppose by contradiction that $b)$ does not hold. Then, using a Lions type result for Orlicz-Sobolev spaces developed  in \cite{AGJ}, we know that
\begin{eqnarray}\label{limB1}
\int_{\mathbb{R}^{N}}B(\vert u_{n}\vert)dx \rightarrow 0.
\end{eqnarray}
By \ref{f1}-\ref{f2}, \ref{R3}-\ref{R4}, \ref{H2} and \ref{B3}, for each $\eta >0$, there exists $c_{\eta} > 0$, such that
\begin{eqnarray}\label{limB2}
\vert f(u_{n})u_{n}\vert \leq \eta c_{1}\Phi(\vert u_{n}\vert) + c_{\eta}b_{2}B(\vert u_{n}\vert).
\end{eqnarray}
Then, gathering the boundedness of  $(u_{n})$ in $L^{\Phi}(\mathbb{R}^{N})$, \eqref{limB1} and \eqref{limB2}, we get
\begin{eqnarray*}
\int_{\mathbb{R}^{N}}f(u_{n})u_{n}dx \rightarrow 0,
\end{eqnarray*}
which leads to
\[
\int_{\mathbb{R}^{N}}\phi(\vert \nabla u_{n}\vert)\vert\nabla u_{n}\vert^{2}dx + \mu\int_{\mathbb{R}^{N}}\phi(\vert u_{n}\vert)\vert u_{n}\vert^{2}dx\rightarrow  0.
\]
The last limit combined with \ref{H2} gives
\[
\int_{\mathbb{R}^{N}}\Phi(\vert \nabla u_{n}\vert)dx + \mu\int_{\mathbb{R}^{N}}\Phi(\vert u_{n}\vert)dx\rightarrow  0,
\]
implying that $u_{n}\rightarrow 0 $ in $Y_{\mu}$.
\end{pf}
\begin{thm} \label{TGS}
Suppose that \ref{H1}-\ref{H3}, \ref{R1}-\ref{R4}, \ref{B1}-\ref{B4}, \ref{f1}-\ref{f3} and $(V_{0})$ hold. Then, \eqref{PA} has a nonnegative  ground state solution.
\end{thm}
\begin{pf}
The result follows combining the mountain pass theorem with Lemmas \ref{la} and \ref{Lemacompaut}.
\end{pf}
\subsection{Behavior asymptotic, regularity and positiveness of solutions}
In this section, we will study the behavior asymptotic, regularity and positiveness of the nonnegatives solutions of \eqref{PA}. Hereafter, we adapt some arguments found in \cite{AF, Fusco, ZF}.
\begin{lem}\label{Lemalimit}
Let $u \in W^{1, \Phi}(\mathbb{R}^{N})$ be a solution of \eqref{PA}, $x_0 \in \mathbb{R}^{N}$ and $R_0>0$. Then,
\[
\int_{A_{k, t}}\vert \nabla u \vert^{\gamma}dx \leq c \Bigg( \int_{A_{k, s}}\Big\vert \frac{u - k}{s- t}\Big\vert^{\gamma^{*}}dx + (k^{\gamma^{*}} + 1)\vert A_{k, s}\vert\Bigg),
\]
where $0< t < s < R_{0}$, $k\geq k_{0}\geq 1$, $k_{0}$ to be chosen and \linebreak $A_{k, \rho}= \{x \in B_{\rho}(x_0) \ : \ u(x) > k \}$.
\end{lem}
\begin{pf}
Let $u$ be a weak solution of \eqref{PA}, $x_{0} \in \mathbb{R}^{N}$ and $R_{0}>1$. Moreover, fix $0< t< s <1<R_{0}$ and $\xi \in C_0^{\infty}(\mathbb{R}^{N})$ verifying
\[
0 \leq \xi\leq 1, \quad supp \xi \subset B_{s}(x_{0}), \quad \xi \equiv1 \ \mbox{on} \ B_{t}(x_{0}) \quad \mbox{and} \quad \vert \nabla \xi\vert\leq \frac{2}{s - t}.
\]
For $k\geq 1$, set $\eta = \xi^{m}(u-k)_{+}$ and
\[
J= \int_{A_{k, s}}\Phi(\vert \nabla u\vert)\xi^{m}dx.
\]
Using $\eta$ as a test function and \ref{H2}, we find
\begin{eqnarray*}
lJ&\leq& m\int_{A_{k, s}}\phi(\vert \nabla u \vert)\vert\nabla u \vert \vert\nabla\xi\vert\xi^{m - 1}(u-k)_{+} dx\\
&&- \mu\int_{A_{k, s}}\phi(\vert u \vert)u\xi^{m}(u-k)_{+} dx + \int_{A_{k, s}}f(u)\xi^{m}(u-k)_{+} dx.
\end{eqnarray*}
Now, from \ref{f1}, \ref{R3}-\ref{R4} and \ref{H2}, there exists a constant $c_{1}>0$ such that
\begin{eqnarray*}
lJ&\leq& m\int_{A_{k, s}}\phi(\vert \nabla u \vert)\vert\nabla u \vert \vert\nabla\xi\vert\xi^{m - 1}(u-k)_{+}dx + c_{1}\int_{A_{k ,s}}b(\vert u\vert)\vert u\vert\xi^{m}(u-k)_{+}dx.
\end{eqnarray*}
For each $\tau \in (0, 1)$, the Young's inequalities gives
\begin{eqnarray*}
\phi(\vert \nabla u \vert)\vert\nabla u \vert \vert\nabla\xi\vert\xi^{m - 1}(u-k)_{+} &\leq& \widetilde{\Phi}(\phi(\vert \nabla u\vert)\vert \nabla u\vert \xi^{m - 1}\tau) + \Phi(\frac{\vert \nabla\xi\vert}{\tau}(u-k)_{+}).
\end{eqnarray*}
It follows from \cite[Lemmas 2.1 and 2.5]{FN}
\[
\Phi\Big(\frac{(u-k)_{+}}{s-t}\Big)\leq \Big\vert \frac{u-k}{s-t}\Big\vert^{m} + 1
\]
and
\[
\widetilde{\Phi}(\phi(\vert \nabla u\vert)\vert \nabla u\vert \xi^{m - 1}\tau)\leq (\xi^{m - 1}\tau)^{\frac{m}{m-1}} \widetilde{\Phi}(\phi(\vert \nabla u\vert)\vert \nabla u\vert).
\]
Now, fixing $\vartheta = \phi(\vert \nabla u \vert)\vert\nabla u \vert \vert\nabla\xi\vert\xi^{m - 1}(u-k)_{+}$, and using the above inequalities  together with the $\Delta_{2}$-condition, we see that
\begin{eqnarray*}
\int_{A_{k, s}}\vartheta dx &\leq&
c_{2}\tau^{\frac{m}{m-1}}\int_{A_{k, s}}\xi^{m}\Phi(\vert \nabla u\vert)dx + c_{3}\int_{A_{k, s}}\Big\vert \frac{u-k}{s-t}\Big\vert^{\gamma^{*}}dx + c_{6}\vert A_{k, s}\vert.
\end{eqnarray*}
Choosing $\tau$ such that $mc_{2}\tau^{\frac{m}{m-1}} < l$, we derive
\begin{eqnarray}\label{eq2P}
J &\leq& c_{1}\int_{A_{k ,s}}b(\vert u\vert)\vert u\vert\xi^{m}(u-k)_{+}dx + c_{3}\int_{A_{k, s}}\Big\vert \frac{u-k}{s-t}\Big\vert^{\gamma^{*}}dx + c_{7}\vert A_{k, s}\vert.\quad \quad
\end{eqnarray}
By Young's inequalities,
\begin{eqnarray*}
b(\vert u\vert)\vert u\vert\xi^{m}(u-k)_{+} &\leq& \widetilde{B}(b(\vert u\vert)\vert u\vert \xi^{m}) + B(\vert u -k\vert),
\end{eqnarray*}
from where it follows that
\[
\int_{A_{k, s}}b(\vert u\vert)\vert u\vert\xi^{m}(u-k)_{+}dx \leq c_{4}\int_{A_{k, s}}B(\vert u -k\vert)dx +  \int_{A_{k, s}}B(\vert k\vert)dx.
\]
Applying \cite[Lemma 2.1]{FN} for function $B$, we find
\[
B\Big(\frac{(u-k)_{+}}{s-t}\Big)\leq \Big\vert \frac{u-k}{s-t}\Big\vert^{b_{2}} + 1.
\]
Now, using that $\gamma^{*}> b_{2}$, we get
\begin{eqnarray*}
\int_{A_{k, s}}b(\vert u\vert)\vert u\vert\xi^{m}(u-k)_{+}dx &\leq& c_{5}\int_{A_{k, s}}\Big\vert \frac{u-k}{s-t}\Big\vert^{\gamma^{*}} dx + c_{8}(k^{\gamma^{*}} + 1)\vert A_{k, s}\vert.\quad
\end{eqnarray*}
From \eqref{eq2P} and the inequality above,
\begin{eqnarray*}
J\leq C\Bigg( \int_{A_{k, s}}\Big\vert \frac{u-k}{s-t}\Big\vert^{\gamma^{*}} dx + (k^{\gamma^{*}} + 1)\vert A_{k, s}\vert\Bigg).
\end{eqnarray*}
If $\Phi \in \mathcal{C}_m$, as $\xi\equiv 1$ on $B_{t}(x_{0})$, we obtain
\[
\int_{A_{k, t}}\vert \nabla u \vert^{m}dx \leq c \Bigg( \int_{{A_{k, s}}}\Big\vert \frac{u - k}{s- t}\Big\vert^{\gamma^{*}}dx + (k^{\gamma^{*}} + 1)\vert A_{k, s}\vert\Bigg).
\]
On the other hand, if $\Phi \notin \mathcal{C}_m$, the inequality
$$
t^{l} \leq \Phi(t)+1 \,\,\, \forall t \geq 0,
$$
leads to
\[
\int_{A_{k, t}}\vert \nabla u \vert^{l}dx \leq c \Bigg( \int_{{A_{k, s}}}\Big\vert \frac{u - k}{s- t}\Big\vert^{\gamma^{*}}dx + (k^{\gamma^{*}} + 1)\vert A_{k, s}\vert\Bigg),
\]
finishing the proof. \end{pf}

\vspace{0.5 cm}

The lemma below is crucial in our work, because it shows that the solutions of \eqref{PA} are in $L^{\infty}(\mathbb{R}^N)$ and they go to zero at infinity.

\begin{lem}\label{L1}
Let $u \in W^{1, \Phi}(\mathbb{R}^{N})$ be a nonnegative solution of \eqref{PA}. Then, $u \in L^{\infty}(\mathbb{R}^{N})$ and $u \in C^{1, \alpha}_{loc}(\mathbb{R}^{N})$. Furthermore,
\[
\lim_{|x| \to +\infty}u(x) =0.
\]
\end{lem}
\begin{pf}
We will divide our proof into three steps.
\begin{flushleft}
\textbf{Step 1: $ u \in L^{\infty}(\mathbb{R}^{N})$}.
\end{flushleft}
To begin with, fix $R_1>0$, $K_{0}\geq 1$ and $x_{0} \in \mathbb{R}^{N}$. Consider $K>K_{0}$ and define
\[
\sigma_{n} = \frac{R_1}{2} + \frac{R_1}{2^{n+1}}, \quad \overline{\sigma}_{n} = \frac{\sigma_{n} + \sigma_{n+1}}{2} \quad \mbox{and}\quad K_{n} = \frac{K}{2}\Big(1-\frac{1}{2^{n+1}} \Big).
\]
For each $n \in \mathbb{N}$, fix
\[
J_{n} = \int_{A_{K_{n}, \sigma_{n}}}(( u-K_{n})_{+})^{\gamma^{*}}dx
\]
and
\[
\xi_{n}= \xi\bigg(\frac{2^{n+1}}{R_1}\Big(\vert x\vert - \frac{R_1}{2}\Big) \bigg) \quad x \in \mathbb{R}^{N},
\]
where $\xi \in C^{1}(\mathbb{R})$ satisfies
\[
0\leq \xi \leq 1,\quad \xi(t) = 1 \ \mbox{for}\ t\leq \frac{1}{2} \quad \mbox{and} \quad \xi(t) = 0 \ \mbox{for} \ t\geq \frac{3}{4}.
\]
From definition of $\xi_{n}$,
\[
\xi_{n} = 1\ \mbox{in} \ B_{\sigma_{n + 1}}(x_{0}) \quad \mbox{and} \quad \xi_{n} = 0 \ \mbox{outside} \ B_{\overline{\sigma}_{n}}(x_{0}),
\]
consequently
\begin{eqnarray*}
J_{n+1} &\leq& 
\int_{B_{R_{1}}(x_{0})}((u - K_{n+1})_{+}\xi_{n})^{\gamma^{*}}dx.
\end{eqnarray*}
As $(u - K_{n+1})_{+}\xi_{n} \in W_{0}^{1, \gamma}(B_{R_{1}}(x_{0}))$, it follows from continuous imbedding and Poincar\'e's inequality that
\[
J_{n+1}^{\frac{\gamma}{\gamma^{*}}} \leq C(N, \gamma)\int_{B_{R_{1}}(x_{0})}\vert \nabla((u - K_{n+1})_{+}\xi_{n} )\vert^{\gamma}dx.
\]
Now, proceeding in the same way as in \cite[Lemma 5.5]{AF}, we find the inequality
\[
J_{n} \leq CA^{\eta}J_{n}^{1 + \eta},
\]
where $C = C(N, \gamma, \gamma^{*}, R_1)$, $A = 2^{(\gamma + \gamma^{*})\frac{\gamma^{*}}{\gamma}}$ and $\eta = \frac{\gamma^{*}}{\gamma}-1$.

We claim that
\begin{eqnarray}\label{J1}
J_{0} \leq C^{\frac{1}{\eta}}A^{-\frac{1}{\eta^{2}}}.
\end{eqnarray}
Indeed, note that,
\[
J_{0} = \int_{A_{K_{0}, \sigma_{0}}}(u - K_{0})_{+}^{\gamma^{*}}dx  \leq \int_{\mathbb{R}^{N}}(u - K_{0})_{+}^{\gamma^{*}}dx.
\]
Then, by the Lebesgue's Theorem, $\displaystyle \lim_{K_{0}\rightarrow +\infty}J_{0}=0$, from where it follows that \eqref{J1} holds for $K_0>0$ large enough. Thus, by \cite[Lemma 4.7]{LU},
\[
\lim_{n \to +\infty}J_{n}= 0.
\]
On the other hand,
\[
\lim_{n \to +\infty}J_{n} =\lim_{n \to +\infty}\int_{A_{K_{n}, \sigma_{n}}}(( u-K_{n})_{+})^{\gamma^{*}}dx=\int_{A_{\frac{K}{2}, \frac{R_1}{2}}}(( u-\frac{K}{2})_{+})^{\gamma^{*}}dx.
\]
Hence,
\[
\int_{A_{\frac{K}{2},  \frac{R_1}{2}}}(( u-\frac{K}{2})_{+})^{\gamma^{*}}dx = 0,
\]
leading to
\[
u(x)\leq \frac{K}{2}\quad \mbox{a.e in}\  B_{\frac{R_1}{2}}(x_{0}).
\]
Recalling that $x_0$ is arbitrary, the last inequality ensures that
\begin{eqnarray*}
u(x)\leq \frac{K}{2}\quad \mbox{a.e in}\  \mathbb{R}^{N},
\end{eqnarray*}
showing that $u \in L^{\infty}(\mathbb{R}^{N})$.

\begin{flushleft}
\textbf{Step 2: $u \in C^{1, \alpha}_{loc}(\mathbb{R}^{N})$.}
\end{flushleft}
This regularity  follows applying the results found in DiBenedetto \cite{Db} and Lieberman \cite{Lb1}.

\begin{flushleft}
\textbf{Step 3: $\displaystyle\lim_{\vert x\vert\rightarrow +\infty}u(x) = 0$.}
\end{flushleft}
Repeating the same arguments used in Step 1 with $K=\frac{\delta}{2}$ and $K_{0} = \displaystyle\frac{\delta}{4}$ for $\delta > 0$, we find
\[
J_{0} = \int_{A_{K_{0}, \sigma_{0}}}(u - K_{0})_{+}^{\gamma^{*}}dx \leq \int_{B_{\sigma_{0}}(x_{0})}(u - \frac{\delta}{4})_{+}^{\gamma^{*}}dx.
\]
Since
\[
\lim_{\vert x_{0} \vert\rightarrow +\infty}\int_{B_{\sigma_{0}}(x_{0})}(u - \frac{\delta}{4})_{+}^{\gamma^{*}}dx = 0,
\]
we also have
\[
\lim_{\vert x_{0} \vert\rightarrow +\infty}J_{0} = 0.
\]
Then, there exists $R_{*}>0$, such that
$$
J_{0}\leq C^{\frac{1}{\eta}}A^{-\frac{1}{\eta^{2}}} \,\,\, \mbox{for} \,\,\, \vert x_{0} \vert > R_{*}.
$$
Applying the Lemma \cite[Lemma 4.7]{LU},
\[
\lim_{n \to +\infty}J_{n}= 0 \quad \mbox{if} \quad \vert x_{0} \vert > R_{*},
\]
that is,
\[
\int_{A_{\frac{\delta}{4},  \frac{R_1}{2}}}(u - \frac{\delta}{4})_{+}^{\gamma^{*}}dx=0\quad \mbox{for} \quad \vert x_{0} \vert > R_{*}.
\]
The continuity of $u$ yields
\[
u(x)\leq \frac{\delta}{4} \quad \forall x \in \mathbb{R}^{N} \setminus B_{R_*}(0),
\]
finishing the proof of the lemma.
\end{pf}

Thanks to Lemma \ref{L1}, we are able to prove the positiveness of the ground state solution obtained in Theorem \ref{TGS}.
\begin{corol}\label{CRI}
Let $u \in Y_{\mu}\backslash \{0\}$ be a nonnegative solution of \eqref{PA}. Then, $u$ is positive solution.
\end{corol}
\begin{pf}
By Lemma \ref{L1},  $u \in L^{\infty}(\mathbb{R}^N)$. If  $\Omega \subset \mathbb{R}^{N}$ is a bounded domain,  the regularity theory found in \cite{Lb1}  implies that $ u \in C^{1}(\overline{\Omega})$. Using this fact, in the sequel, we fix $M_{1} > \max\{\|\nabla u \|_{\infty},1\}$ and
\begin{align}
\widetilde{\phi}(t)=\left\{
\begin{array}
[c]{rcl}%
\phi(t), \quad t \leq M_{1}\\
\displaystyle\frac{\phi(M_{1})}{M_{1}^{l-2}}t^{l-1}, \quad t\geq M_{1}.
\end{array}\nonumber
\right.
\end{align}
A simple computation yields there exists $\alpha_{1}>0$ verifying
\begin{eqnarray}\label{har1}
\widetilde{\phi}(\vert y\vert)\vert y \vert^{2}\geq \alpha_{1} \vert y \vert^{l}-\alpha_{1} \,\,\, \forall y \in \mathbb{R}^{N}.
\end{eqnarray}
Using the above function, we also define the function
\[
\widetilde{A}(y) = \frac{1}{\alpha_1}\widetilde{\phi}(\vert y\vert)y,
\]
which verifies
$$
\widetilde{A}(y)=\frac{1}{\alpha_1}\phi(|y|)y  \,\,\, \mbox{for} \,\,\, |y| \leq M_1
$$
and
\begin{eqnarray*}\label{har2}
\vert \widetilde{A}(y) \vert \leq c \vert y\vert^{l-1} \,\,\, \forall y \in \Omega,
\end{eqnarray*}
for some positive constant $c$. By \eqref{har1},
\begin{eqnarray*}\label{har3}
y\widetilde{A}(y) =\alpha^{-1}\widetilde{\phi}(\vert y\vert)\vert y\vert^{2} \geq \vert y \vert^{l}- 1 \,\,\, \forall y \in \Omega.
\end{eqnarray*}
Setting $B(x) = \displaystyle\frac{1}{\alpha_1}\big(\mu\phi(\vert u(x)\vert)u(x) - f(u(x))\big)$, we infer that $u$ is a weak solution of the quasilinear problem
\begin{eqnarray*}
-div \widetilde{A}(\nabla u(x)) + B(x) = 0 \,\,\, \mbox{in} \,\,\, \Omega.
\end{eqnarray*}
Since $\Omega$ is arbitrary, by \cite[Theorem 1.1]{Harnack}, we deduce that $u$ is a positive solution.
\end{pf}

\section{Existence of positive ground state solution for the Nonautonomous case}
In this section, our first goal is to prove the existence of positive ground state solution for the following quasilinear problem
\begin{align}
\left\{
\begin{array}
[c]{rcl}%
- \Delta_{\Phi}u + V(\epsilon x)\phi(|u|)u & = & f(u)~  \mbox{in}~ \mathbb{R}^{N},\\
u(x) >0 \,\, \mbox{in} \,\, \mathbb{R}^{N}, & & \\
u \in W^{1, \Phi}(\mathbb{R}^{N}), &  &
\end{array}
\right. \tag{$ P_{\epsilon} $}\label{Pe}%
\end{align}
when $\epsilon$ is small enough. To this end, let us denote by $I_{\epsilon} : X_{\epsilon} \rightarrow \mathbb{R}$ the energy functional related to \eqref{Pe} given by
\[
I_{\epsilon}(u)= \int_{\mathbb{R}^{N}}\Phi(\vert\nabla u\vert)dx + \int_{\mathbb{R}^{N}}V(\epsilon x)\Phi(\vert u\vert)dx - \int_{\mathbb{R}^{N}}F(u)dx,
\]
where $X_{\epsilon}$ denotes the subspace of $W^{1, \Phi}(\mathbb{R}^{N})$ given by
$$
X_{\epsilon}=\left\{ u \in W^{1, \Phi}(\mathbb{R}^{N}) \,:\, \int_{\mathbb{R}^{N}}V(\epsilon x)\Phi({\vert u\vert})dx < +\infty \right\},
$$
endowed with the norm
\[
\Vert u\Vert_{\epsilon} = \Vert \nabla u\Vert_{\Phi} + \Vert u\Vert_{\Phi, V_{\epsilon}},
\]
where
\[
\Vert u\Vert_{\Phi, V_{\epsilon}}:= \inf\Big\{ \lambda > 0;  \int_{\mathbb{R}^{N}}V(\epsilon x)\Phi\Big(\frac{\vert u\vert}{\lambda}\Big)dx \leq 1\Big\}.
\]
By $(V_0)$,  the embeddings
\[
X_{\epsilon}\hookrightarrow L^{\Phi}( \mathbb{R}^{N}) \quad \mbox{and}\quad X_{\epsilon}\hookrightarrow L^{B}( \mathbb{R}^{N})
\]
are continuous. Using the above embeddings, it is standard to show that $I_{\epsilon} \in C^{1}(X_{\epsilon}, \mathbb{R})$ with
\[
I^{'}_{\epsilon}(u)v = \int_{\mathbb{R}^{N}}\phi(\vert\nabla u\vert)\nabla u \nabla vdx + \int_{\mathbb{R}^{N}}V(\epsilon x)\phi(\vert u\vert)uvdx - \int_{\mathbb{R}^{N}}f(u)vdx,
\]
for all $u,v \in X_{\epsilon}$. Thereby, $u \in X_\epsilon$ is a weak solution of (\ref{Pe}) if, and only if, $u$ is a critical point of $I_{\epsilon}$. Furthermore, by (\ref{POSITIVIDADE}),  the critical points of $I_{\epsilon}$ are nonnegative.

Moreover, as in \cite[Lemma 4.1]{AGJ},  $I_{\epsilon}$ verifies the mountain pass geometry. Therefore, applying again a  version of Mountain Pass Theorem without the Palais-Smale condition found in \cite{Willen}, there exists a $(PS )_{c_{\epsilon}}$ sequence $(u_{n}) \subset X_{\epsilon}$, that is, a sequence satisfying
\[
I_{\epsilon}(u_{n})\rightarrow c_{\epsilon} \quad \mbox{and}\quad I_{\epsilon}^{'}(u_{n})\rightarrow 0,
\]
where $c_{\epsilon}$ is the mountain level of $I_{\epsilon}$. Arguing as in Lemma \ref{Lcar}, we can prove that $c_{\epsilon}$ verifies the  equality below
\[
c_{\epsilon} = \inf_{u \in \mathcal{N}_{\epsilon}}I_{\epsilon}(u) = \inf_{u \in X_{\epsilon}\backslash\{0\}}\max_{t\geq 0}I_{\epsilon}(tu),
\]
where $\mathcal{N}_{\epsilon}$ denotes the Nehari manifold related to $I_{\epsilon}$ given by
\[
\mathcal{N}_{\epsilon} = \big\{ u \in X_{\epsilon}\backslash \{0 \} \ : \ I^{'}_{\epsilon}(u)u = 0 \big\}.
\]

Our first lemma shows that the Nehari manifolds has a positive distance from the origin in $X_\epsilon$ and its proof  follows by using standard arguments.
\begin{lem} \label{Nehari}
For all $u \in \mathcal{N}_{\epsilon}$, there exists $k>0$, which is independent of $\epsilon$, such that
\[
\Vert u\Vert_{\epsilon} > k.
\]
\end{lem}

The lemma below establishes some properties involving the $(PS)$ sequences for $I_\epsilon$. Since it follows repeating the same arguments used in Lemma \ref{la}, we will omit its proof.

\begin{lem}\label{Lemaltda}
Let $(u_{n})$ be a sequence $(PS)_{c}$. Then, $(u_{n})$ is a bounded sequence in $X_{\epsilon}$. Moreover, if $u_n \rightharpoonup u$ in $X_\epsilon$, for some subsequence, we have that
\begin{description}
  \item[$(i)$] $ u_{n}(x)\rightarrow  u(x)\quad \mbox{a. e. in} \ \mathbb{R}^{N}$,
  \item[$(ii)$] $ \nabla u_{n}(x)\rightarrow \nabla u(x)\quad \mbox{a. e. in} \ \mathbb{R}^{N}$,
  \item[$(iii)$] $I^{'}_{\epsilon}(u) = 0$,
  \item[$(iv)$] ${u_{n}}_{-} = \max\{-u_{n}, 0\}\rightarrow 0$ in $X_{\epsilon}$.
	\item[$(v)$] ${u_{n}}_{+} = \max\{u_{n}, 0\}$ is a $(PS)_c$ sequence for $I_\epsilon$. Hence, we can assume that $u_n \geq 0$ for all $n \in \mathbb{N}$, and so, $u \geq 0$.
\end{description}
\end{lem}

The next lemma is a version of Lemma \ref{Lemacompaut} with $E_\mu$ replaced by $I_\epsilon$, and the proof can be made using the same ideas.

\begin{lem}\label{Lemacomp}
Let $(u_{n}) \subset X_{\epsilon}$ be a $(PS)_{c}$ sequence for $I_{\epsilon}$ with $u_{n}\rightharpoonup 0 $ in $X_{\epsilon}$. Then,
\begin{description}
  \item[$a)$] $u_{n}\rightarrow 0$ in  $X_{\epsilon}$,
\end{description}	
	or
\begin{description}
  \item[$b)$] There exist constants $R_0, \tau>0$ such that
  \begin{eqnarray*}
\liminf_{n \rightarrow +\infty}\sup_{y \in \mathbb{R}^{N}}\int_{B_{R_0}(y)}\Phi(\vert u_{n}\vert)dx \geq \tau >0.
\end{eqnarray*}
\end{description}
\end{lem}

The next lemma is crucial to show that $I_\epsilon$ verifies the $(PS)$ condition at some levels.

\begin{lem}\label{lema1}
Assume that $V_{\infty} < +\infty$ and let $(v_{n})\subset X_{\epsilon}$ be a $(PS)_{c}$ sequence for $I_{\epsilon}$ with $v_{n}\rightharpoonup 0$ in $X_{\epsilon}$. If $v_{n}\nrightarrow 0$ in $X_{\epsilon}$, then $d_{V_{\infty}}\leq c$.
\end{lem}
\begin{pf}
Let $(t_{n}) \subset (0, +\infty)$ be a sequence such that $t_{n}v_{n} \in \mathcal{M}_{V_{\infty}}$. Adapting the same arguments explored in \cite[Lemma 3.3]{AG1}, we conclude that $(t_{n})$ satisfies
\[
\limsup_{n \rightarrow +\infty}t_{n}\leq 1.
\]
Next, we will study the following cases:
\begin{flushleft}
\textbf{Case 1.} \, $\displaystyle\limsup_{n \rightarrow +\infty} t_{n} = 1$.
\end{flushleft}
In this case, there exists a subsequence of $(t_{n})$, still denoted by itself, such that $t_{n}\rightarrow 1$. Using the fact that $(v_{n})$ is bounded in $X_{\epsilon}$ together with the condition $(V_0)$, given $\eta>0$, there is $R_0>0$ such that
\begin{eqnarray*}
I_{\epsilon}(v_{n}) - I_{V_{\infty}}(t_{n}v_{n}) &\geq& \int_{\mathbb{R}^{N}}\Big[\Phi(\vert \nabla v_{n}\vert)-\Phi(\vert\nabla t_{n}v_{n}\vert)\Big]dx + \int_{\mathbb{R}^{N}}\Big[F(t_{n}v_{n}) - F(v_{n})\Big]dx\\
&& + V_{\infty}\int_{B^{c}_{R_0}(0)}\Big[ \Phi(\vert v_{n}\vert) - \Phi(\vert t_{n}v_{n}\vert)\Big]dx -\eta C_4 +o_{n}(1).
\end{eqnarray*}
From the mean value theorem,
\[
\int_{\mathbb{R}^{N}}\Big[\Phi(\vert\nabla v_{n}\vert)-\Phi(\vert \nabla t_{n}v_{n}\vert)\Big]dx = o_{n}(1), \quad \int_{\mathbb{R}^{N}}\Big[\Phi(\vert v_{n}\vert)-\Phi(\vert t_{n}v_{n}\vert)\Big]dx = o_{n}(1)
\]
\and
\[
\int_{\mathbb{R}^{N}}\Big[F(t_{n}v_{n}) - F(v_{n})\Big]dx = o_{n}(1).
\]
Then,
\[
c =I_{\epsilon}(v_{n})+ o_{n}(1) \geq d_{V_{\infty}} + I_{\epsilon}(v_{n}) - I_{V_{\infty}}(t_{n}v_{n}) + o_{n}(1) \geq d_{V_{\infty}} - C_4\eta + o_{n}(1).
\]

As $\eta$ is arbitrary, taking the limit of $n \to +\infty$, we derive that $c \geq d_{V_{\infty}}$.

\begin{flushleft}
\textbf{Case 2.} \, $\displaystyle\limsup_{n \rightarrow +\infty} t_{n} = t_{0}<1$.
\end{flushleft}
In this case, we may suppose that there exists a subsequence of $(t_{n})$, still denoted by $(t_{n})$, satisfying $t_{n}\rightarrow t_{0}$ and $t_{n}<1$ for all $n \in \mathbb{N}$. By \ref{H3} and \ref{f3}, the functions $P(s)=\Phi(s) - \frac{1}{m}\phi(s)s^{2}$  and $Q(s)=\frac{1}{m}f(s)s - F(s)$ are increasing in $(0,+\infty)$, thus
$$
d_{V_{\infty}}\leq\int_{\mathbb{R}^{N}}P(\vert\nabla v_{n}\vert)dx + \int_{\mathbb{R}^{N}}V(\epsilon x)P(\vert v_{n}\vert)dx
+2\eta\int_{\mathbb{R}^{N}}\Phi(\vert v_{n}\vert)dx + \int_{\mathbb{R}^{N}}Q(v_{n})dx,
$$
implying  that
\[
d_{V_{\infty}} \leq c_{1}\eta + I_{\epsilon}(v_{n}) - \frac{1}{m}I^{'}_{\epsilon}(v_{n})v_{n}=  C_{4}\eta + c + o_{n}(1).
\]
As $\eta$ is arbitrary, taking the limit of $n \to +\infty$, it follows that $c \geq d_{V_{\infty}}$.
\end{pf}

The next two results are technical and they will use to show the local Palais-Smale condition for $I_\epsilon$. The first one is an immediate consequence of a result Brezis-Lieb in \cite{BL} and it has the ensuing statement
\begin{lem}
Assume that $\phi$ satisfies \ref{H1}-\ref{H3}. Suppose also that $(u_{n}) \subset X_{\epsilon}$ is a $(PS)_{c}$ sequence for $I_{\epsilon}$ with
$$
u_{n}(x) \to u(x) \,\,\, \mbox{and} \,\,\, \nabla u_n(x) \to \nabla u(x) \,\,\, \mbox{a.e in} \,\,\, \mathbb{R}^{N}.
$$
Then
\begin{eqnarray}\label{ig1}
\int_{\mathbb{R}^{N}} \Phi(\vert\nabla v_{n}\vert)dx = \int_{\mathbb{R}^{N}} \Phi(\vert\nabla u_{n}\vert)dx - \int_{\mathbb{R}^{N}} \Phi(\vert\nabla u\vert)dx  + o_{n}(1),\quad
\end{eqnarray}
\begin{eqnarray}\label{ig2}
\int_{\mathbb{R}^{N}}V(\epsilon x) \Phi(\vert v_{n}\vert)dx = \int_{\mathbb{R}^{N}} V(\epsilon x)\Phi(\vert u_{n}\vert)dx - \int_{\mathbb{R}^{N}} V(\epsilon x)\Phi(\vert u\vert)dx  + o_{n}(1),\quad
\end{eqnarray}
and
\begin{eqnarray}\label{ig3}
\int_{\mathbb{R}^{N}}F(v_{n})dx = \int_{\mathbb{R}^{N}} F(u_{n})dx - \int_{\mathbb{R}^{N}} F(u)dx  + o_{n}(1),\quad
\end{eqnarray}
where $v_{n} = u_{n} - u$.
\end{lem}
\begin{prop}\label{Pess}
Let $(u_{n}) \subset X_{\epsilon}$ be a $(PS)_{c}$ sequence for $I_{\epsilon}$ such that $u_{n}\rightharpoonup u$ in $X_{\epsilon}$. Then
\begin{eqnarray}\label{propa}
I_{\epsilon}(v_{n}) = I_{\epsilon}(u_{n}) - I_{\epsilon}(u) + o_{n}(1)
\end{eqnarray}
and
\begin{eqnarray}\label{propb}
I^{'}_{\epsilon}(v_{n}) = o_{n}(1),
\end{eqnarray}
where $v_{n} = u_{n} - u$.
\end{prop}
\begin{pf}
The proof of \eqref{propa} follows from \eqref{ig1}, \eqref{ig2} and \eqref{ig3}. To prove \eqref{propb}, it is sufficient to show that
\[
\Vert I^{'}_{\epsilon}(v_{n}) - I^{'}_{\epsilon}(u_{n}) + I^{'}_{\epsilon}(u)\Vert = o_{n}(1).
\]
Let $w \in X_{\epsilon}$ with $\Vert w\Vert_{\epsilon} \leq 1$ and  $Q_{\epsilon}(w) =\vert \big(I^{'}_{\epsilon}(v_{n}) - I^{'}_{\epsilon}(u_{n}) + I^{'}_{\epsilon}(u)\big)w\vert $. Note that
\begin{eqnarray*}
Q_{\epsilon}(w) &\leq& \int_{\mathbb{R}^{N}}\big\vert\phi(\vert\nabla (v_{n} + u)\vert)\nabla (v_{n} + u) - \phi(\vert\nabla v_{n}\vert)\nabla v_{n} -\phi(\vert\nabla u\vert)\nabla u \big\vert \vert\nabla w\vert dx\\
&&+\int_{\mathbb{R}^{N}}\big\vert\phi(\vert v_{n} + u\vert)(v_{n} + u) - \phi(\vert v_{n}\vert) v_{n} -\phi(\vert u\vert) u \big\vert \vert w \vert V(\epsilon x)dx\\
&& + \int_{\mathbb{R}^{N}}\big\vert (f(v_{n} + u) - f(v_{n}) -f(u))w\big\vert dx.
\end{eqnarray*}
Applying the Proposition \ref{Teoape1}, we obtain
\begin{eqnarray}\label{y1}
Q_{\epsilon}(w) &\leq& o_{n}(1) \Vert  w\Vert_{\epsilon} + \int_{\mathbb{R}^{N}}\big\vert f(u_{n})w - f(v_{n})w -f(u)w\big\vert dx.
\end{eqnarray}
Since $w \in X_{\epsilon}$, using \ref{f1}, \ref{R3} and H\"{o}lder inequality combined with compact embedding and Lebesgue's Theorem yields, for each $\eta>0$, there exists $R_0>0$ such that 
$$
\int_{\mathbb{R}^{N}}\big\vert f(u_{n})w - f(v_{n})w -f(u)w\big\vert dx \leq o_{n}(1) + c_{1}\eta +\int_{B^{c}_{R_0}(0)}\big\vert f(u_{n})w - f(v_{n})w\big\vert dx.
$$
Now, it is enough to show that
\begin{eqnarray}\label{y2}
\int_{B^{c}_{R_0}(0)}\big\vert f(v_{n} + u)w - f(v_{n})w\big\vert dx = o_{n}(1),\, \, \mbox{uniformly for}\, \, \Vert w \Vert_{\epsilon}\leq 1.
\end{eqnarray}
By \ref{f1}, \ref{R1}-\ref{R2} and \ref{B1}-\ref{B2},
$$
\vert  f(v_{n} + u) - f(v_{n}) \leq c_{1}r(\vert v_{n}\vert )\vert u\vert + c_{1}r(\vert u \vert)\vert u \vert + c_{2}b(\vert v_{n}\vert)\vert u\vert + c_2b(\vert u\vert)\vert u\vert.
$$
Now, using \ref{R4} and \ref{B4}, we derive that
\[
\sup_{\Vert w\Vert_{\epsilon}\leq 1}\int_{B^{c}_{R_0}(0)}\big\vert f(v_{n} + u)w - f(v_{n})w\big\vert dx = o_{n}(1).
\]
Hence, by \eqref{y1} and \eqref{y2},
\[
\sup_{\Vert w\Vert_{\epsilon}\leq 1}Q_{\epsilon}(w) = o_{n}(1),
\]
implying that
\[
\Vert I_{\epsilon}^{'}(v_{n})\Vert = o_{n}(1).
\]
\end{pf}

The next result establishes a compact embedding involving the space $X_{\epsilon}$. The particular cases where $\Phi(t)=\frac{1}{p}|t|^{p}$ for $N>p>1$ and $\Phi(t)=\frac{1}{p}|t|^{p}+\frac{1}{q}|t|^{q}$ for $N>q>p>1$ were considered in \cite{David} and \cite{AG1} respectively.
\begin{lem}\label{Lemaic}
If $V_{\infty} = +\infty$, the embedding $X_{\epsilon}\hookrightarrow L^{\Phi}(\mathbb{R}^{N})$ is compact. In particular, the embedding $X_{\epsilon}\hookrightarrow L^{B}(\mathbb{R}^{N})$ is compact, if $B$ is a N-function satisfying \ref{B4}.
\end{lem}
\begin{pf}
Let $(u_{n}) \subset X_{\epsilon}$ be a sequence bounded with $u_n \rightharpoonup u$ in $X_\epsilon$. For each $\eta >0$, there exists $R_0>0$ such that
\[
V(\epsilon x)> \frac{1}{\eta}, \quad |x|> R_0.
\]
Thereby,
$$
\int_{\mathbb{R}^{N}}\Phi(\vert u_{n} - u\vert )dx \leq  \eta\int_{B^{c}_{R_0}(0)}V(\epsilon x)\Phi(\vert u_{n} - u\vert)dx + o_{n}(1).
$$
Here, we have used that
$$
u_{n}\rightarrow u \quad \mbox{in}\quad L^{\Phi}(B_{R_0}(0)).
$$
Using the boundedness of $(u_n)$ and the arbitrariness of $\eta$, we can ensure that
$$
u_{n}\rightarrow u \,\,\, \mbox{in} \,\,\, L^{\Phi}(\mathbb{R}^{N}).
$$
To conclude our proof, if $B$ verifies \ref{B4}, for each $\eta>0$, there is $C_{\eta}>0$ such that
$$
\int_{\mathbb{R}^{N}}B(|u_n-u|)dx \leq C_\eta \int_{\mathbb{R}^{N}} \Phi(|u_n-u|)dx +\eta \int_{\mathbb{R}^{N}}\Phi_*(|u_n-u|)dx, \,\,\, \forall n \in \mathbb{N}.
$$
As $\eta$ is arbitrary, using the boundedness of $(u_n-u)$ in $L^{\Phi_*}(\mathbb{R}^{N})$ together with the fact that $u_n \to u$ in
$L^{\Phi}(\mathbb{R}^{N})$, we deduce that
$$
\int_{\mathbb{R}^{N}}B(|u_n-u|)dx \to 0,
$$
showing that $u_n \to u$ in $L^{B}(\mathbb{R}^{N})$. \end{pf}

The proposition below shows where the $(PS)$ condition holds for $I_\epsilon$.

\begin{prop}\label{prop1}
Assume that $V_{\infty}< +\infty$. Then, $I_{\epsilon}$ satisfies the $(PS)_{c}$ condition for any $c < d_{V_{\infty}}$. If $V_{\infty} = +\infty$, the $(PS)_{c}$ condition occurs  for any $c \in \mathbb{R}$.
\end{prop}
\begin{pf}
Let $(u_{n}) \subset X_{\epsilon}$ be a sequence such that $I_{\epsilon}(u_{n})\rightarrow c$ and $I^{'}_{\epsilon}(u_{n})\rightarrow 0$. Setting $v_{n} = u_{n} - u$, it follows from  Lemma  \ref{Lemaltda} and  Proposition \ref{Pess}
$$
I_{\epsilon}(v_{n})= c- I_{\epsilon}(u) + o_{n}(1)= \overline{c} + o_{n}(1) \quad \mbox{and} \quad I^{'}_{\epsilon}(v_{n}) = o_{n}(1).
$$
By \ref{H2} and \ref{f2},
\begin{eqnarray*}
I_{\epsilon}(u)&=& I_{\epsilon}(u) - \frac{1}{m}I^{'}(u)u = \int_{\mathbb{R}^{N}}\big(\Phi(\vert \nabla u\vert)-\frac{1}{m}\phi(\vert \nabla u\vert)\vert\nabla u\vert^{2}\big)dx\\
&&+ \int_{\mathbb{R}^{N}}\big(\frac{1}{m}f(u)u - F(u)\big)dx +\int_{\mathbb{R}^{N}}V(x)\big(\Phi(\vert u\vert)-\frac{1}{m}\phi(\vert u\vert)\vert u\vert^{2}\big)dx\geq 0.
\end{eqnarray*}
Assume that $V_{\infty} < +\infty$. As $(v_{n})$ is a sequence $(PS)_{\overline{c}}$ and $\overline{c} \leq c< d_{V_{\infty}}$, the Lemma \ref{lema1} yields
\[
v_{n}\rightarrow 0 \quad \mbox{in} \  X_{\epsilon},
\]
showing that $u_n \to u$ in $X_\epsilon$. Now, if $V_{\infty} = +\infty$, by Lemma \ref{Lemaic}, up to a subsequence, $v_{n}\rightarrow 0$ in the spaces $L^{\Phi}(\mathbb{R}^{N})$ and $L^{B}(\mathbb{R}^{N})$. By \ref{f1}, \ref{R3}-\ref{R4}, \ref{H2} and \ref{B3}, we derive that
\begin{eqnarray*}
\int_{\mathbb{R}^{N}}f(v_{n})v_{n}dx \rightarrow 0.
\end{eqnarray*}
The boundedness of $(v_n)$ combined with the limit $I^{'}_{\epsilon}(v_{n}) = o_n(1)$ gives
\begin{eqnarray}\label{qc}
\int_{\mathbb{R}^{N}}\phi(\vert \nabla v_{n}\vert )\vert\nabla v_{n}\vert^{2}dx + \int_{\mathbb{R}^{N}}V(x)\phi(\vert v_{n}\vert)\vert v_{n}\vert^{2}dx = o_{n}(1).
\end{eqnarray}
Gathering (\ref{qc}) with $(\phi_2)$, we conclude that $v_n\rightarrow 0$ in $X_{\epsilon}$,
or equivalently, $ u_{n}\rightarrow u \quad \mbox{in} \  X_{\epsilon}.$ \end{pf}

The next proposition shows that the $(PS)$ condition holds for $I_\epsilon$ on $\mathcal{N}_\epsilon$ at some levels.

\begin{prop}\label{prop2}
Assume that $V_{\infty} < + \infty$. Then, $I_{\epsilon}$ satisfies $(PS)_{c}$ condition on $\mathcal{N}_{\epsilon}$ for $c \in (0, d_{V_{\infty}})$. If $V_{\infty} = + \infty$, the $(PS)$ condition occurs for any $c \in \mathbb{R}$.
\end{prop}
\begin{pf}
Let $(u_{n})\subset \mathcal{N}_{\epsilon}$ be a sequence verifying
\[
I_{\epsilon}(u_{n})\rightarrow c \quad \mbox{and}\quad \parallel I^{'}_{\epsilon}(u_{n})\parallel_{*} = o_{n}(1).
\]
Then, there exists $(\lambda_{n})\subset \mathbb{R}$ such that
\begin{eqnarray}\label{ML}
I^{'}_{\epsilon}(u_{n}) = \lambda_{n}J^{'}_{\epsilon}(u_{n}) + o_{n}(1),
\end{eqnarray}
where $J_{\epsilon}: X_{\epsilon}\rightarrow \mathbb{R}$ is given by
\[
J_{\epsilon}(u) = \int_{\mathbb{R}^{N}}\phi(\vert \nabla u\vert)\vert \nabla u\vert^{2}dx + \int_{\mathbb{R}^{N}}V_{\epsilon}(x)\phi(\vert u\vert)\vert u\vert^{2}dx - \int_{\mathbb{R}^{N}}f(u)udx.
\]
From definition of $J_{\epsilon}$,
\begin{eqnarray*}
J^{'}_{\epsilon}(u_{n})u_{n} &=& \int_{\mathbb{R}^{N}}\big(\phi^{'}(\vert\nabla u_{n}\vert)\vert\nabla u_{n}\vert + 2 \phi(\vert\nabla u_{n}\vert)\big)\vert\nabla u_{n}\vert^{2}dx\\
&&+ \int_{\mathbb{R}^{N}}V_{\epsilon}(x)\big(\phi^{'}(\vert u_{n}\vert)\vert u_{n}\vert + 2 \phi(\vert u_{n}\vert)\big)\vert u_{n}\vert^{2}dx\\
&& - \int_{\mathbb{R}^{N}}\big( f^{'}(u_{n})(u_{n})^{2} + f(u_{n})u_{n}\big)dx.
\end{eqnarray*}
The last equality, \ref{H3} and $I_{\epsilon}^{'}(u_{n})u_{n} = 0$ combine to give
\begin{eqnarray*}
-J^{'}_{\epsilon}(u_{n})u_{n} &\geq& \int_{\mathbb{R}^{N}}\big(f^{'}(u_{n})(u_{n})^{2}-(m-1)f(u_{n})u_{n}\big)dx.
\end{eqnarray*}
Since  $(u_{n})$ is a bounded sequence, we can assume that $u_{n}\rightharpoonup u$ in $X_{\epsilon}$. Moreover, there exists a subset $\Omega\subset \mathbb{R}^{N}$ with positive measure such that $u > 0$ a.e. in $\Omega$. Suppose that
\[
\limsup_{n\rightarrow +\infty}J^{'}_{\epsilon}(u_{n})u_{n} = 0.
\]
Using \ref{f3} and applying the Fatous' Lemma in the last inequality, we get
$$
0 \geq \int_{\Omega}\big(f^{'}(u)(u)^{2}-(m-1)f(u)u\big)dx>0,
$$
which is an absurd. Hence $\displaystyle\limsup_{n\rightarrow +\infty}J^{'}_{\epsilon}(u_{n})u_{n}< 0$, implying that $\lambda_{n} = o_{n}(1)$ for some subsequence. This combined with \eqref{ML} gives $I^{'}_{\epsilon}(u_{n}) = o_{n}(1)$, \linebreak  showing that $(u_{n})$ is $(PS)_{c}$. Now, the result follows from  Proposition \ref{prop1}.
\end{pf}

As a by product of the arguments used in the proof of the last proposition, we have the following result.

\begin{corol}\label{cor1}
The critical points of functional $I_{\epsilon}$ on $\mathcal{N}_{\epsilon}$ are critical points of $I_{\epsilon}$ in $X_{\epsilon}$.
\end{corol}

Using the previous result proved in this section, we are able to prove the existence of positive ground state solution for  \eqref{Pe} when $\epsilon$ small enough.
\begin{thm}
Assume that \ref{H1}-\ref{H5}, \ref{R1}-\ref{R4}, \ref{B1}-\ref{B4}, \ref{f1}-\ref{f3} and $(V_0)$ hold. Then, there exists $\overline{\epsilon} > 0$ such that \eqref{Pe} has a nonnegative ground state solution $u_{\epsilon}$ for all $0 < \epsilon <\overline{\epsilon}$.
\end{thm}
\begin{pf}
We begin recalling that $I_{\epsilon}$ satisfies the mountain pass geometry. Then, there exists $(u_{n}) \subset X_{\epsilon}$ satisfying
$$
I_{\epsilon}(u_{n})\rightarrow c_{\epsilon} \,\,\, \mbox{and} \,\,\, I^{'}_{\epsilon}(u_{n})\rightarrow 0,
$$
where $c_\epsilon$ denotes the mountain pass level associated with $I_\epsilon$.
If $V_{\infty} = +\infty$, by Proposition \ref{prop1}, there exists $u \in X_{\epsilon}$ such that, for some subsequence, $
u_{n}\rightarrow u \quad \mbox{in} \quad X_\epsilon.$ Therefore,
\[
I_{\epsilon}(u) = c_{\epsilon} \quad \mbox{and} \quad I^{'}_{\epsilon}(u) = 0.
\]

If $V_{\infty} < +\infty$, we must prove that $c_\epsilon \in (0, d_{V_{\infty}})$ for $\epsilon$ small enough. To do that, consider without loss of generality that $V(0)=V_{0}$ and $\mu \in (V_{0}, V_{\infty})$. Then, by definition of the mountain pass level, we have that
\begin{eqnarray*}
d_{V_{0}}< d_{\mu} < d_{V_{\infty}}.
\end{eqnarray*}
Thus, there exists $w \in Y_{\mu}$ with compact support, such that
$$
E_{\mu}(w) = \displaystyle\max_{t \geq 0}E_{\mu}(tw)  \,\,\, \mbox{and} \,\,\, E_{\mu}(w) < d_{V_{\infty}}.
$$
From $(V_0)$, there exists $\overline{\epsilon}>0$ such that
\[
V(\epsilon x) <\mu, \quad \forall \epsilon \in (0, \overline{\epsilon})\quad \mbox{and}\quad \forall x \in supp w.
\]
By a direct computation, $\displaystyle\max_{t>0}I_{\epsilon}(tw) \leq E_{\mu}(w)$, and so,
\[
c_{\epsilon} <  d_{V_{\infty}}, \,\,\,\,\,\, \forall \epsilon \in (0, \overline{\epsilon}).
\]
Now, the theorem follows from Proposition \ref{prop1}.
\end{pf}
\section{Multiplicity of solutions to \eqref{Pe}}
In this section, our main goal is to show the existence of multiple solutions, and the behavior of their maximum points of the solutions in relation to the set $M$.

\subsection{Preliminary results}
In order to prove our main theorem, we must fix some notations. In what follows, $w\in W^{1, \Phi}(\mathbb{R}^{N})$ denotes a positive ground state solution of $(P_{V_{0}})$, $\delta>0$ is small enough and $\eta \in C_{0}^{\infty}([0, +\infty),[0,1])$ is given by
\[
\eta(s) =
\left\{
\begin{array}{rcl}
  1, &&\mbox{if} \quad 0 \leq s \leq \frac{\delta}{2}\\
  0, &&\mbox{if} \quad s \geq \delta.\\
\end{array}
\right.
\]
Using the above notations, for each $y \in M$,  we set
\[
\Psi_{\epsilon, y}(x)=\eta(\vert\epsilon x -y\vert)w(\frac{\epsilon x -y}{\epsilon}).
\]
Note that, $\Psi_{\epsilon, y} \in X_{\epsilon}$ for each $y \in M$. Thus, there exists $t_{\epsilon}>0$ such that $t_{\epsilon}\Psi_{\epsilon, y} \in {\mathcal N}_{\epsilon}$ and
\[
I_{\epsilon}(t_{\epsilon}\Psi_{\epsilon, y})= \max_{t\geq 0}I_{\epsilon}(t\Psi_{\epsilon, y}).
\]
Consequently, we can define $\widetilde{\Psi}_{\epsilon}: M \rightarrow {\mathcal N}_{\epsilon}$ by $\widetilde{\Psi}_{\epsilon}(y) = t_{\epsilon}\Psi_{\epsilon, y}$. By construction, $\widetilde{\Psi}_{\epsilon}$ has compact support for any $y \in M$.
\begin{lem}\label{lemapsi}
The function $\widetilde{\Psi}_{\epsilon}$ verifies the following limit
\[
\lim_{\epsilon\rightarrow 0}I_{\epsilon}(\widetilde{\Psi}_{\epsilon}(y)) = d_{V_{0}}, \quad \mbox{uniformly in} \quad y \in M.
\]
\end{lem}
\begin{pf}
It is sufficient to show that for each $(y_{n}) \subset M$ and $(\epsilon_{n})\subset \mathbb{R^{+}}$ with $\epsilon_{n} \rightarrow 0$, there is a subsequence such that
\begin{eqnarray*}
I_{\epsilon}(\widetilde{\Psi}_{\epsilon_{n}}(y_{n})) \rightarrow d_{V_{0}}.
\end{eqnarray*}
Recall firstly that $J^{'}_{\epsilon_{n}}(\widetilde{\Psi}_{\epsilon_{n}}(y_{n}))\widetilde{\Psi}_{\epsilon_{n}}(y_{n})= 0 $, that is,
\begin{eqnarray*}
\int_{\mathbb{R}^{N}}\widehat{\phi}(\vert\nabla(\widetilde{\Psi}_{\epsilon_{n}}(y_{n})\vert)dx + \int_{\mathbb{R}^{N}}V(\epsilon_{n}x)\widehat{\phi}(\vert \widetilde{\Psi}_{\epsilon_{n}}(y_{n})\vert)dx =  \int_{\mathbb{R}^{N}}f(\widetilde{\Psi}_{\epsilon_{n}}(y_{n}))\widetilde{\Psi}_{\epsilon_{n}}(y_{n})dx,
\end{eqnarray*}
where $\widehat{\phi}(s)=\phi(s)s^{2}$  for all $s \geq 0$. Using \ref{H2} and \cite[Lemma 2.1]{FN},
$$
\begin{array}{l}
\displaystyle \int_{\mathbb{R}^{N}}\widehat{\phi}(\vert\nabla(\widetilde{\Psi}_{\epsilon_{n}}(y_{n})\vert)dx +\int_{\mathbb{R}^{N}}V(\epsilon_{n}z + y_{n})\widehat{\phi}(\vert \widetilde{\Psi}_{\epsilon_{n}}(y_{n})\vert)dx\\
\leq m \xi_{1}(t_{\epsilon_{n}})\Big[\displaystyle \int_{\mathbb{R}^{N}}\Phi(\vert\nabla(\Psi_{\epsilon_{n}, y_{n}})\vert)dx + \int_{\mathbb{R}^{N}}V(\epsilon_{n}z + y_{n})\Phi(\vert\Psi_{\epsilon_{n}, y_{n}}\vert)dx\Big],
\end{array}
$$
where $\xi_{1}(t)=\max\{t^{l}, t^{m}\}$. On the other hand, considering the change of variable $z = \displaystyle\frac{\epsilon_{n}x - y_n}{\epsilon_n}$, it follows that
\[
\int_{\mathbb{R}^{N}}f(\widetilde{\Psi}_{\epsilon_{n}}(y_{n}))\widetilde{\Psi}_{\epsilon_{n}}(y_{n})dx =\int_{\mathbb{R}^{N}}f(t_{\epsilon_{n}}\eta(|\epsilon_{n}z|)w(z))t_{\epsilon_{n}}\eta(|\epsilon_{n}z|)w(z)dx.
\]
As  $\eta\equiv1$ on $B_{\frac{\delta}{2}}(0)$ and $B_{\frac{\delta}{2}}(0) \subset B_{\frac{\delta}{2\epsilon_{n}}}(0)$, we have that
\begin{eqnarray*}
\int_{B_{\frac{\delta}{2}}(0)}\frac{f(t_{\epsilon_{n}}w(z))}{(t_{\epsilon_{n}}w(z))^{m-1}}\vert t_{\epsilon_{n}}w(z)\vert^{m}dx &\leq& m\xi_{1}(t_{\epsilon_{n}})\Big[\int_{\mathbb{R}^{N}}\Phi(\vert\nabla(\Psi_{\epsilon_{n}, y_{n}})\vert)dx\\
&&+ \int_{\mathbb{R}^{N}}V(\epsilon_{n}z + y_{n})\Phi(\vert\Psi_{\epsilon_{n}, y_{n}}\vert)dx\Big].
\end{eqnarray*}
By Lemma \ref{L1}, we know that $w$ is a continuous function. Then, there is $z_{0} \in \mathbb{R}^{N}$ such that
\[
w(z_{0})=\min_{z\in B_{\frac{\delta}{2}}(0)}w(z),
\]
and so,
\begin{eqnarray*}
\frac{f(t_{\epsilon_{n}}w(z_{0}))}{(t_{\epsilon_{n}}w(z_{0}))^{m-1}}\int_{B_{\frac{\delta}{2}}(0)}\vert t_{\epsilon_{n}}w(z)\vert^{m}dx &\leq& m\xi_{1}(t_{\epsilon_{n}})\Big[\int_{\mathbb{R}^{N}}\Phi(\vert\nabla(\Psi_{\epsilon_{n}, y_{n}})\vert)dx\\
&& +\int_{\mathbb{R}^{N}}V(\epsilon_{n}z + y_{n})\Phi(\vert\Psi_{\epsilon_{n}, y_{n}}\vert)dx\Big].
\end{eqnarray*}
By \ref{f3}, there are $d_{1}, d_{2}>0$ such that
\begin{eqnarray*}
&&\big[d_{1}(t_{\epsilon_{n}}w(z_{0}))^{\theta-m} - d_{2}(t_{\epsilon_{n}}w(z_{0}))^{-m}\big]t_{\epsilon_{n}}^{m}\int_{B_{\frac{\delta}{2}}(0)}\vert w(z)\vert^{m}dx\\
&&\leq m\xi_{1}(t_{\epsilon_{n}})\Big[\int_{\mathbb{R}^{N}}\Phi(\vert\nabla(\Psi_{\epsilon_{n}, y_{n}})\vert)dx + \int_{\mathbb{R}^{N}}V(\epsilon_{n}z + y_{n})\Phi(\vert\Psi_{\epsilon_{n}, y_{n}}\vert)dx\Big].
\end{eqnarray*}
Now, arguing by contradiction, we will suppose that, for some subsequence,
\[
t_{\epsilon_{n}}\rightarrow +\infty \,\,\, \mbox{and} \,\,\, t_{\epsilon_{n}} \geq 1 \,\,\, \forall n \in \mathbb{N}.
\]
Thereby, $\xi_{1}(t_{\epsilon_{n}})=t_{\epsilon_{n}}^{m}$ and
\begin{eqnarray*}
&&\big[d_{1}(t_{\epsilon_{n}}w(z_{0}))^{\theta-m} - d_{2}(t_{\epsilon_{n}}w(z_{0}))^{-m}\big]\int_{B_{\frac{\delta}{2}}(0)}\vert w(z)\vert^{m}dx\\
&&\leq m\Big[\int_{\mathbb{R}^{N}}\Phi(\vert\nabla(\Psi_{\epsilon_{n}, y_{n}})\vert)dx + \int_{\mathbb{R}^{N}}V(\epsilon_{n}z + y_{n})\Phi(\vert\Psi_{\epsilon_{n}, y_{n}}\vert)dx\Big].
\end{eqnarray*}
Note that the right side of the above inequality is bounded, while that
$$
[d_{1}(t_{\epsilon_{n}}w(z_{0}))^{\theta-m} - d_{2}(t_{\epsilon_{n}}w(z_{0}))^{-m}\big] \to +\infty,
$$
what is a contradiction. Therefore $(t_{\epsilon_{n}})$ is bounded, and for some subsequence, there exists $t_{0}\geq 0$ such that
$t_{\epsilon_{n}} \rightarrow t_{0}.$ Now, recalling that $\widetilde{\Psi}_{\epsilon_{n}}(y_{n}) \in \mathcal{N}_{\epsilon_n}$, we know that $\|\widetilde{\Psi}_{\epsilon_{n}}(y_{n})\|_{\epsilon} \geq k$ for all $\epsilon>0$. Using again the Lebesgue's Theorem, it is possible to prove that
$$
E'_{V_0}(t_0w)(t_0w)=0 \,\,\, \mbox{and} \,\,\,\,  \|t_0w\|_{Y_{V_0}}\geq  k,
$$
implying that $t_0>0$ and $t_0w \in \mathcal{M}_{V_0}$. However, as $w$ is a ground state solution, we must have $t_0=1$. Since $t_n \to 1$,
we apply again the Lebesgue's Theorem to get
\[
\displaystyle\lim_{n\rightarrow +\infty}I_{\epsilon}(\widetilde{\Psi}_{\epsilon_{n}}(y_{n})) = E_{V_{0}}(w)= d_{V_{0}},
\]
finishing the proof.
\end{pf}

In the sequel, we state an important result of compactness for for $E_\mu$ restricts to Nehari Manifolds $\mathcal{M}_\mu$. The proof is quite similar to the case where $\Phi(t)=\frac{1}{p}|t|^{p}$, which was made in \cite{AG2006}, then we will omit its proof.
\begin{lem}\label{LC}\textbf{(A Compactness Lemma)}
Let $(u_{n}) \subset \mathcal{M}_{\mu}$ be a sequence satisfying $E_{\mu}(u_{n})\rightarrow d_{\mu}$. Then,
\begin{description}
  \item[$a)$] $(u_{n})$ has a subsequence strongly convergent in $Y_{\mu}$, or
  \item[$b)$] there exists a sequence $(\widetilde{y}_{n}) \subset \mathbb{R}^{N}$ such that, up to a subsequence, $v_{n}(x)=u_{n}(x+\widetilde{y}_{n})$ converges strongly in $Y_{\mu}$.
\end{description}
In particular, there exists a minimizer for $d_{\mu}$.
\end{lem}

Our next result will be crucial to study the concentration of the solutions.

\begin{prop}\label{Propimp}
Let $\epsilon_{n}\rightarrow 0$ and $(u_{n})\subset {\mathcal N}_{\epsilon_{n}}$ be such that $I_{\epsilon_{n}}(u_{n})\rightarrow d_{V_{0}}$. Then, there exists a sequence $(\widetilde{y}_{n})\subset\mathbb{R}^{N}$, such that $v_{n}(x) = u_{n}(x+\widetilde{y}_{n})$ has a convergent subsequence in $Y_{V_{0}}$. Moreover, up to a subsequence, $y_{n}\rightarrow y \in M$, where $y_{n}=\epsilon_{n}\widetilde{y}_{n}$.
\end{prop}
\begin{pf}
Applying the Lemma \ref{Lemacomp}, we find a sequence $(\widetilde{y}_{n}) \subset \mathbb{R}^{N}$ and positive constants $R_0$ and $\tau$ such that
\begin{eqnarray}\label{lz}
\int_{B_{R_0}(\widetilde{y}_{n})}\Phi(u_{n})dx\geq \tau >0.
\end{eqnarray}
Since, $(u_{n})\subset \mathcal{N}_{\epsilon_{n}}$ and $I_{\epsilon_{n}}(u_{n})\rightarrow d_{V_{0}}$, we know that $(u_{n})$ is bounded in $Y_{V_{0}}$. Then, setting  $v_{n}(x) = u_{n}(x +\widetilde{y}_{n})$, up to a subsequence, there is $v \in Y_{V_{0}}\setminus\{0\}$ such that $ v_{n}\rightharpoonup v \quad \mbox{in} \ Y_{V_{0}}$  and  $v \not=0$ by (\ref{lz}).  Let $t_{n} > 0$ such that $\widetilde{v}_{n} = t_{n}v_{n} \in \mathcal{M}_{V_{0}}$. Then,
$$
E_{V_{0}}(\widetilde{v}_{n}) \leq I_{\epsilon}(t_{n}u_{n}) \leq \max_{t\geq 0}I_{\epsilon_{n}}(tu_{n}) = I_{\epsilon_{n}}(u_{n}),
$$
and so,
\begin{eqnarray}\label{limitvn}
E_{V_{0}}(\widetilde{v}_{n})\rightarrow d_{V_{0}} \quad \mbox{and} \quad (\widetilde{v}_{n}) \subset \mathcal{M}_{V_{0}}.
\end{eqnarray}
Using \eqref{limitvn}, it follows $\widetilde{v}_{n} \rightharpoonup \widetilde{v}$ in $Y_{V_{0}}$. It is possible to check that $(t_{n})$ is bounded. Moreover, up to a subsequence, $t_{n} \rightarrow t_{0}>0$, implying that $\widetilde{v} = t_{0}v$. Thus, it follows that $\widetilde{v}_{n} \rightarrow \widetilde{v}$ in $Y_{V_{0}}$ showing that $v_{n}\rightarrow v$ in $Y_{V_{0}}$. Now, the proof follows as in \cite[Proposition 4.1]{AG2006}.
\end{pf}

In the sequel, for any $\delta > 0$, let $\rho = \rho(\delta)>0$ be such that $M_{\delta} \subset B_{\rho}(0)$. Let $\chi: \mathbb{R}^{N} \rightarrow \mathbb{R}^{N}$ given by
\[
\chi(x) =
\left\{
\begin{array}{rcl}
  x, \quad \mbox{if} \quad x \in B_{\rho}(0)\\
  \displaystyle\frac{\rho x}{\vert x\vert}, \quad \mbox{if} \quad x \in B^{c}_{\rho}(0)\\
\end{array}
\right.
\]
and $\beta: \mathcal{N}_{\epsilon} \rightarrow \mathbb{R}^{N}$ given by
$$
\beta(u)= \frac{\displaystyle\int_{\mathbb{R}^{N}}\chi(\epsilon x)\Phi(\vert u\vert)dx}{\displaystyle\int_{\mathbb{R}^{N}}\Phi(\vert u\vert)dx}.
$$

\begin{lem}
The function $\widetilde{\Psi}_{\epsilon_{n}}$ satisfies the following limit
\[
\lim_{\epsilon_{n}\rightarrow 0}\beta(\widetilde{\Psi}_{\epsilon_{n}}(y)) = y, \ \mbox{uniformemente em}\ M.
\]
\end{lem}
\begin{pf}
The lemma follows by using the definition of $\widetilde{\Psi}_{\epsilon_{n}}(y)$ together with the Lebesgue's Theorem.
\end{pf}

In the sequel, we consider the function $h: \mathbb{R}^{+} \rightarrow \mathbb{R}^{+}$ given by
\[
h(\epsilon)=\sup_{y \in M}|I_\epsilon(\widetilde{\Psi}_{\epsilon}(y))-d_{V_0}|,
\]
which verifies, by Lemma \ref{lemapsi}, $\displaystyle \lim_{\epsilon \to 0}h(\epsilon)=0$. Moreover, we set
\[
\widetilde{\mathcal{N}}_{\epsilon}:= \{u \in \mathcal{N}_{\epsilon} \ : \ I_{\epsilon}(u) \leq d_{V_{0}} + h(\epsilon) \}.
\]
From Lemma \ref{lemapsi}, $\widetilde{\Psi}_{\epsilon}(y) \in \widetilde{\mathcal{N}}_{\epsilon}$, showing that $\widetilde{\mathcal{N}}_{\epsilon}\neq\emptyset$. Using the above notation, we have the following result:
\begin{lem}\label{lemapsi1}
Let $\delta>0$ and $M_{\delta}= \{ x \in \mathbb{R}^{N}\ :\ dist(x, M)\leq \delta\}$. Then, the limit below hols
\[
\lim_{\epsilon\rightarrow 0}\sup_{u \in \widetilde{\mathcal{N}}_{\epsilon}}\inf_{y \in M_{\delta}} \big\vert\beta(u) - y \big\vert = 0.
\]
\end{lem}
\begin{pf}
The proof follows as in \cite[Lemma 4.4]{AG2006}.
\end{pf}
\section{Proof of Theorem~\ref{T1} }
We will divide the proof into two parts:
\subsection{Part I: Multiplicity of solutions}
For each $\epsilon > 0$ small enough, the Lemmas \ref{lemapsi} and \ref{lemapsi1} ensure that $\beta\circ\widetilde{\Psi}_{\epsilon}$ is homotopic to inclusion map $id : M \rightarrow M_{\delta}$, showing the following estimate holds
\[
cat_{\widetilde{\mathcal{N}}_{\epsilon}}(\widetilde{\mathcal{N}}_{\epsilon})\geq cat_{M_{\delta}}(M).
\]
Since  $I_{\epsilon}$ satisfies the $(PS)_{c}$ condition for $c \in (d_{V_{0}}, d_{V_{0}} + h(\epsilon))$, by the Lusternik-Schnirelman theory of critical points, we can conclude that $I_{\epsilon}$ has at least $cat_{M_{\delta}}(M)$ critical points on $\mathcal{N}_{\epsilon}$. Consequently by Corollary \ref{cor1}, $I_{\epsilon}$ has at least $cat_{M_{\delta}}(M)$ critical points in $X_{\epsilon}$.
\subsection{Part II: The behavior of maximum points}
The next result plays an important role in the study of the behavior of the maximum points of the
solutions.
\begin{prop}\label{pd}
Let $v_{j} \in W^{1, \Phi}(\mathbb{R}^{N})$ be a solution of the following problem
\begin{align}
\left\{
\begin{array}
[c]{rcl}%
- \Delta_{\Phi}v_{j} + V_{j}(x)\phi(\vert v_{j}\vert)v_{j} & = & f(v_{j})~  \mbox{in}~ \mathbb{R}^{N},\\
v_{j} \in W^{1, \Phi}(\mathbb{R}^{N}), &  &\\
v_{j} > 0 \quad in \quad \mathbb{R}^{N},& &
\end{array}
\right. \tag{$ P_{j} $}\label{Pj}%
\end{align}
where $V_{j}(x) = V(\epsilon_{j}x + \epsilon_{j}\widetilde{y}_{j})$,  $\epsilon_{j}\widetilde{y}_{j} \to y$ in $\mathbb{R}^{N}$ and $v_j \to v$ in $W^{1,\Phi}(\mathbb{R}^{N})$. Then $v_{j} \in C^{1, \gamma}_{loc}(\mathbb{R}^{N})$, $v_{j}\in L^{\infty}(\mathbb{R}^{N})$ and there exists $C > 0$ such that $ \|v_{j}\|_\infty \leq C$ for all $j \in \mathbb{N}$. Furthermore
\[
\lim_{\vert x \vert \rightarrow +\infty} v_{j}(x) = 0, \quad \mbox{uniformly in j}.
\]
\end{prop}

In order to prove the above proposition, we will use the lemma below, which is a version of Lemma \ref{Lemalimit} for the nonautonomous case, whose the proof explores the same ideas.

\begin{lem}
If $v_{j} \in W^{1, \Phi}(\mathbb{R}^{N})$ is a solution of problem \eqref{Pj}, $x_0 \in \mathbb{R}^{N}$ and $R_0>0$, then
\[
\int_{A_{j, k, t}}\vert \nabla v_{j} \vert^{\gamma}dx \leq c \Bigg( \int_{A_{j,k, s}}\Big\vert \frac{v_{j} - k}{s- t}\Big\vert^{\gamma^*}dx + (k^{\gamma^*} + 1)\vert A_{j, k, s}\vert\Bigg),
\]
where $0< t < s < R_{0}$, $k\geq k_{0}\geq 1$, $k_{0}$ to be chosen, and \linebreak $A_{j,k, \rho}= \{x \in B_{\rho}(x_0) \ : \ v_j(x) > k \}$.
\end{lem}

\begin{flushleft}
\textbf{Proof of Proposition \ref{pd}}
\end{flushleft}
\begin{pf}
To begin with, fix  $\sigma_{n}, \overline{\sigma}_{n}$ and $K_{n}$ as in Lemma \ref{L1}.
\begin{flushleft}
\textbf{Step 1:} There is $C>0$ such that $\|v_j\|_{\infty} \leq C \,\,\, \forall j \in \mathbb{N}.$
\end{flushleft}
Define for each $j$
\[
J_{n, j} = \int_{A_{j,K_{n}, \sigma_{n}}}(( v_{j}-K_{n})_{+})^{\gamma^*}dx.
\]
Repeating the same arguments found in proof of Lemma \ref{L1}, we see that for each $j \in\mathbb{N}$,
\[
J_{n, j} \leq CA^{\eta}J_{n, j}^{1 + \eta} \,\,\,\,\,\, \forall n \in \mathbb{N},
\]
where $C$ is independent of $n$. Now, we claim that
\[
J_{0, j} \leq C^{\frac{1}{\eta}}A^{-\frac{1}{\eta^{2}}},\quad \mbox{for} \quad j\approx +\infty.
\]
Indeed, since $v_{j} \rightarrow v$ in $W^{1, \Phi}(\mathbb{R}^{N})$,
\[
\limsup_{K_{0}\rightarrow +\infty}\left(\limsup_{j\rightarrow +\infty}J_{0, j}\right) =\limsup_{K_{0}\rightarrow +\infty}\left(\limsup_{j\rightarrow +\infty} \int_{A_{j,K_{0}}, \sigma_{0}}(v_j - K_{0})_{+}^{\gamma^*}dx \right)=0.
\]
Then, there are $j_0 \in \mathbb{N}$ and $K^*_0 >0$ such that
\[
J_{0, j} \leq C^{\frac{1}{\eta}}A^{-\frac{1}{\eta^{2}}},\quad \mbox{for} \quad j\geq j_0 \,\,\, \mbox{and} \,\,\, K_0 \geq K^*_0.
\]
By Lemma \cite[Lemma 4.7]{LU},
\[
\lim_{n \to +\infty}J_{n, j}=0, \quad  \mbox{for} \quad j \geq j_0.
\]
On the other hand,
\[
\lim_{n \to +\infty}J_{n, j} = \lim_{n \to +\infty}\int_{A_{j,K_{n}, \sigma_{n}}}(( v_{j}-K_{n})_{+})^{\gamma^*}dx =\int_{A_{j,\frac{K}{2}, \frac{R_1}{2}}}(( v_{j}-\frac{K}{2})_{+})^{\gamma^*}dx.
\]
Thus,
\[
\int_{A_{j,\frac{K}{2},  \frac{R_1}{2}}}(( v_{j}-\frac{K}{2})_{+})^{\gamma^*}dx = 0,\quad \mbox{for all} \quad j \geq j_0 ,
\]
leading to
\[
v_{j}(x)\leq \frac{K}{2}\quad \mbox{a.e in} \,\,\, B_{\frac{R_1}{2}}(x_{0}),\quad \mbox{for all} \quad j \geq j_0.
\]
Since $x_0 \in \mathbb{R}^{N}$ is arbitrary, we deduce that
\begin{eqnarray*}
v_{j}(x)\leq \frac{K}{2}\quad  \mbox{a.e in} \,\,\, \mathbb{R}^{N}, \quad \mbox{for all} \quad j \geq j_0,
\end{eqnarray*}
that is,
$$
\|v_j\|_\infty \leq \frac{K}{2} \,\,\, \forall j \geq j_0.
$$
Setting $C=\max\{\frac{K}{2},\|v_1\|_\infty, .....,\|v_{j_0-1}\|_\infty \}$, we have that
$$
\|v_j\|_\infty \leq C \,\,\,\,\,\,\,  \forall j \in \mathbb{N}.
$$

\begin{flushleft}
\textbf{Step 2: $v_j \in C^{1, \alpha}_{loc}(\mathbb{R}^{N})$.}
\end{flushleft}
This regularity  follows applying results found in DiBenedetto \cite{Db} and Lieberman \cite{Lb1}.

\begin{flushleft}
\textbf{Step 3:} $\displaystyle\lim_{\vert x\vert\rightarrow +\infty}v_{j}(x) = 0$ uniformly in $j$.
\end{flushleft}
Repeating the same arguments explored in Step 1, for each $\delta > 0$ we have that
\[
\limsup_{|x_0|\to +\infty}\left(\limsup_{j \to +\infty}J_{0, j}\right) =\limsup_{|x_0|\to +\infty}\left(\limsup_{j \to +\infty}J_{0, j} \int_{A_{j,K_{0}, \sigma_{0}}}(v_{j} - \frac{\delta}{4})_{+}^{\gamma^*}dx \right)=0.
\]
Thereby, there are $R_{*}>0$ and $j_0 \in \mathbb{N}$ such that
$$
J_{0, j}\leq C^{\frac{1}{\eta}}B^{-\frac{1}{\eta^{2}}} \,\,\,\, \mbox{when} \,\,\,\, \vert x_{0} \vert > R_{*} \,\,\, \mbox{for} \,\,\, j \geq j_0.
$$
Applying the Lemma \cite[Lemma 4.7]{LU},
\[
\lim_{n \to +\infty}J_{n, j}=0 \quad \mbox{if} \quad \vert x_{0} \vert > R_{*}, \quad \mbox{for} \quad j \geq j_0 ,
\]
that is,
\[
\int_{B_{\frac{R_1}{2}}(x_{0})}(v_{j} - \frac{\delta}{4})_{+}^{\gamma^*}dx=0\quad \mbox{for} \quad \vert x_{0} \vert > R_{*}, \quad \mbox{for all} \quad j \geq j_0
\]
showing that
\[
v_{j}(x)\leq \frac{\delta}{4} \quad \mbox{for} \ \ x \in B_{\frac{R_1}{2}}(x_{0}) \quad \mbox{and} \quad \vert x_{0} \vert > R_{*}, \quad \mbox{for all} \quad j \geq j_0.
\]
Now, increasing $R_*$ if necessary,
\[
v_{j}(x)\leq \frac{\delta}{4} \quad \mbox{for} \quad \vert x \vert > R_{*} \quad \mbox{for all} \quad j \in \mathbb{N},
\]
finishing the proof. \end{pf}

\begin{corol}
Let $v_{j} \in W^{1, \Phi}(\mathbb{R}^{N})\backslash \{0\}$ a solution of \eqref{Pj}. Then, $v_{j}$ is positive.
\end{corol}
\begin{pf}
The proof follows by using similar arguments as in Corollary \ref{CRI}.
\end{pf}
\begin{lem}\label{Le2}
There exists $\alpha >0$ such that $ \|v_{j}\|_\infty>\alpha$ for all $j \in \mathbb{N}$.
\end{lem}
\begin{pf}
The proof follows arguing by contradiction and using the fact that $\|v_j\|_{\epsilon_j} \geq k$ for all $j \in \mathbb{N}$, where $k$ was given in Lemma \ref{Nehari}. 
\end{pf}

\noindent \textbf{The behavior of maximum points} \, Finally, as in \cite{AG2006}, if $u_{\epsilon_{n}}$ is a solution of problem $(P_{\epsilon_{n}})$. Then, $v_{n}(x) = u_{\epsilon_{n}}(x + \widetilde{y}_{n})$ is a solution of problem
\begin{align}
\left\{
\begin{array}
[c]{rcl}%
- \Delta_{\Phi}v_{n} + V_{n}(x)\phi(\vert v_{n}\vert)v_{n} & = & f(v_{n})~  \mbox{in}~ \mathbb{R}^{N},\\
v_{n} \in W^{1, \Phi}(\mathbb{R}^{N}), &  &\\
v_{n} > 0 \quad in \quad \mathbb{R}^{N}, & &
\end{array}
\right. \tag{$ P_{n} $}\label{P3}%
\end{align}
where $V_{n}(x) = V(\epsilon_{n}x + \epsilon_{n}\widetilde{y}_{n})$ and $\widetilde{y}_{n}$ were given in Proposition \ref{Propimp}. Moreover, up to a subsequence, $v_{n} \rightarrow v$ in $W^{1, \Phi}(\mathbb{R}^{N})$ and $y_{n} \rightarrow y$ in $M$, where $y_{n}=\epsilon_{n}\widetilde{y}_{n}$. Applying the Proposition \ref{pd} and Lemma \ref{Le2}, there exists $q_{n} \in B_{R_{0}}(0)$ such that $v_{n}(q_{n}) = \displaystyle\max_{z \in \mathbb{R}^{N}}v_{n}(z)$, for some $R_0>0$. Hence, $x_{n} = q_{n} + \widetilde{y}_{n}$ is a maximum point of $u_{\epsilon_{n}}$ and
\[
\epsilon_{n}x_{n} \rightarrow y.
\]
Since V is a continuous function, it follows that
\[
\lim_{n \to +\infty}V(\epsilon_{n}x_{n})=V(y) = V_{0},
\]
showing the concentration of the maximum points of the solutions near to minimum points of $V$. \fim

\end{document}